\documentclass{article}
\usepackage{amssymb}
\usepackage{amssymb}
\usepackage{amsmath}
\usepackage[dvips]{graphicx}
\usepackage{amsmath,amssymb}
\usepackage{color,eurosym}

\usepackage{times}
\usepackage[T1]{fontenc}
\newtheorem{teo}{Theorem}[section]

\newtheorem{lem}{Lemma}[section]
\newtheorem{cor}{Corollary}[section]

\newcommand{\be}{\begin{equation}}
\newcommand{\ee}{\end{equation}}
\newcommand{\ba}{\begin{array}}
\newcommand{\ea}{\end{array}}
\newcommand{\bee}{\begin{eqnarray*}}
\newcommand{\eee}{\end{eqnarray*}}
\newcommand{\bea}{\begin{eqnarray}}
\newcommand{\eea}{\end{eqnarray}}

\newcommand{\ef}{e_{{\cal N}}}
\newcommand{\efk}{e_{{\cal N}_{k}}}
\newcommand{\eg}{\tilde{e}_{{\cal N}}}
\newcommand{\egk}{\tilde{e}_{{\cal N}_{k}}}
\newcommand{\fnk}{f_{{\cal N}_{k}}}
\newcounter{algo}[section]
\renewcommand{\thealgo}{\thesection.\arabic{algo}}
\newcommand{\algo}[3]{\refstepcounter{algo}
\begin{center}\begin{figure}[htb]
\framebox[\textwidth]{
\parbox{0.95\textwidth} {\vspace{\topsep}
{\bf Algorithm \thealgo : #2}\label{#1}\\
\vspace*{-\topsep} \mbox{ }\\
{#3} \vspace{\topsep} }}
\end{figure}\end{center}}
\def\IR{\hbox{\rm I\kern-.2em\hbox{\rm R}}}
\def\IC{\hbox{\rm C\kern-.58em{\raise.53ex\hbox{$\scriptscriptstyle|$}}
    \kern-.55em{\raise.53ex\hbox{$\scriptscriptstyle|$}} }}

\begin {document}
\title{Subsampled Inexact Newton methods for minimizing large  sums of convex functions}

\author{Stefania Bellavia
 \thanks{ Department of Industrial Engineering, University of Florence, Viale Morgagni, 40/44, 50134 Florence, Italy, e-mail: {\tt
stefania.bellavia@unifi.it}. Research supported by  Gruppo Nazionale per il Calcolo Scientifico,
(GNCS-INdAM) of Italy. }
 \and Nata\v sa Kreji\' c\thanks{Department of Mathematics and
Informatics, Faculty of Sciences,  University of Novi Sad, Trg
Dositeja Obradovi\'ca 4, 21000 Novi Sad, Serbia, e-mail: {\tt
natasak@uns.ac.rs}. Research supported by Serbian Ministry of Education Science and Technological Development, grant no. 174030.}
 \and Nata\v sa Krklec Jerinki\'c \thanks{Department of Mathematics and
Informatics, Faculty of Sciences,  University of Novi Sad, Trg
Dositeja Obradovi\'ca 4, 21000 Novi Sad, Serbia, e-mail: {\tt
natasa.krklec@dmi.uns.ac.rs}. Research supported by Serbian Ministry of Education, Science and Technological Development, grant no. 174030.}}
\date{November 14, 2018}
\maketitle

\begin{abstract}
This paper deals with the minimization of large sum of convex functions by Inexact Newton (IN) methods employing  subsampled functions, gradients and Hessian approximations. The Conjugate Gradient method is used to compute the inexact Newton step and global convergence is enforced by a nonmonotone line search procedure. The aim is to obtain methods with affordable costs and fast convergence.  Assuming strongly convex functions, R-linear convergence and worst-case iteration complexity of the procedure are investigated when functions and gradients are approximated with increasing accuracy.
A set of rules for the forcing parameters and subsample  Hessian sizes are derived that ensure local q-linear/superlinear convergence of the proposed method.
 The random choice of the Hessian subsample is also considered and  convergence in the mean square, both for finite and infinite sums of functions, is proved. Finally,  global convergence with asymptotic $R$-linear rate of IN methods is extended to the case of sum of convex function and strongly convex objective function.
Numerical results  on well known binary classification problems are also given.  Adaptive strategies for selecting forcing terms and Hessian  subsample size, streaming out of the theoretical analysis, are employed and the numerical results showed that they yield  effective IN methods.
\end{abstract}

\section{Introduction}
The problem we consider is
\be \label{problem1}
\min f_{\cal N}(x) = \frac{1}{N} \sum_{i=1}^{N} f_i(x)
\ee
with $ x \in \mathbb{R}^n, $  $ N$ - very large and all functions $ f_i: {\mathbb R}^n \rightarrow {\mathbb R} $ convex  and $f_{\cal N}$  strongly convex.  We are also interested in the case of large  dimension  $ n. $  There is a number of important problems of this type. To start with, one can be interested in minimizing the objective function stated in the form of mathematical expectation, $ f(x) = E[F(x,w)], $ with $ w $ being a random variable from some probability space.  Given that the analytical expression for mathematical expectation is rarely available, one possibility is to approximate the expectation with Sample Average Approximation, SAA, function. In that case, a sample $ \{w^1,\ldots,w^N\} $ is generated and the approximate objective function of the form (\ref{problem1}) with $ f_i(x) = F(x,w^i) $ is minimized. To ensure good approximation of the original objective function in general, one has to take a very large sample and thus calculating $ f_{\cal N}(x), $ its gradient and Hessian is expensive.

 Binary and multi-class classification problems,  e.g.,  employing softmax activation function and  cross-entropy loss can also be  expressed in the form (\ref{problem1}).  For a given set of data (usually very large), we are interested in classifying the data
 according to a set of rules specified by the data features. This problem has attracted a lot of attention recently due to its importance in machine learning.
As the data set is generally very large it is imperative to use methods that can be implemented with reasonable costs.

In the framework of classical optimization, (\ref{problem1}) is a convex problem that can be solved either by a first order or a second order method. However, the size of $ N $ makes classical approaches prohibitively costly and thus calls for specific methods. One possibility is to consider different subsampling schemes which are used to reduce the cost of calculating $ f_{\cal N}, $ its gradient and Hessian. There are many approaches in the literature, based on the idea of using a small sample subset at the beginning of iterative process and increasing the sample size as the solution is approached, ranging from simple heuristic approach \cite{Friedlander, noc2} to  elaborate schemes that take into account the progress achieved up to the current iteration, \cite{Bastin, B. Tuan, B. teorija, ebnkjm, Feris,nas,nas1,nkjm,  P2, Polak}.

Whatever scheduling one adopts, it is assumed that the full sample is eventually reached, at least for the objective function and thus, the next question  to be discussed is the choice of method. First order methods are attractive due to their low cost. One successful example is the stochastic gradient method, that employs a smaller subset of gradient components and thus reduce the cost even further, \cite{Friedlander}. On the other hand several papers investigate the use of second order methods in this framework and demonstrate advantages in some important problems if the second order methods are correctly implemented, \cite{bgm, bbn, bbn_2017, noc1,noc3,em,mit1,mit2,Pilanci,xu1,xu2}.    For a comprehensive discussion of this issue one can see \cite{nocedalsurvey} and references cited therein.

In this paper we focus on  subsampled Inexact Newton (IN)  methods for (\ref{problem1}) wrapped in a nonmonotone linesearch strategy. As it is well known, in Inexact (or Truncated) Newton methods \cite{DES,nash}, the Newton equation  is approximately solved and in case of large scale problems an iterative  Krylov method is used to compute an approximate solution of the Newton equation. The  convexity of the objective function we are dealing with allows us to use  CG method \cite{kelley}. 

The choice of the nonmonotone strategy is motivated by the fact that  the method uses approximate functions, at least initially, before the full sample is reached. Then, enforcing strict decrease in the Armijo rule might require unnecessary small steps.
We adopt the nonmonotone line search procedure introduced in \cite{LF} and, assuming that each of the functions $ f_i $ is strongly convex, we prove  R-linear convergence. Also the worst-case iteration  complexity is investigated and it is proved that the worst-case complexity bound  of this class of non-monotone algorithms, analyzed in  \cite{grapiglia},
is maintained provided that errors in gradient and function  also decay with R-rate. Namely, the method requires at most
$O(\log(\epsilon^{-1}))$ iterations to reach $f(x^k)-f(x^*)<\epsilon$, where $x^*$ is the minimizer of  \eqref{problem1}.

Then, we turn our attention to the local properties of the method with the aim to obtain a local convergence rate faster than the R-linear convergence provided by  first-order methods.
The local convergence rate of IN methods with full Hessian depends on the choice of forcing terms which governs the error in solving each Newton linear system  \cite{EW1}.
Here, as  the Hessian of the objective function given in  (\ref{problem1}) might be prohibitively expensive to compute, we concentrate on subsampled Hessian and IN method that employs such Hessian approximations.  We  point out that in this context, it  is pointless to solve the Newton equation exactly as the Hessian is generally approximated with a lower accuracy than function and gradient and the Newton model employed is actually a subsampled Newton model. Therefore, the use of CG method that allows us to control the accuracy in the solution of the Newton equation is advisable even if $n$ is not large \cite{noc2}.

Assuming that  the sample size scheduling is given for the objective function and the gradient, i.e.  assuming that eventually one reaches the full sample size $ N $, we analyze  the local convergence of a subsampled Hessian IN method.  The  analysis provides   bounds on the Hessian accuracy requirements that depend on the employed forcing terms.
Adaptive forcing terms, streaming out of the iterative process itself  are derived as well.
Furthermore, it is shown that the local method combines well with the nonmonotone line search, i.e., the  $q$-linear/$q$-superlinear convergence rate  of the local method  are preserved.

In the second part of the paper we consider a randomized method obtained by relaxing the conditions for Hessian subsampling. Hence we prove the q-linear/superlinear convergence in the mean sense assuming that the Hessian approximation is good enough with high probability. The analysis yields relation between the Hessian subsample size, the forcing term  and the (computable) sampled  gradient at each iteration.     Q-linear convergence in the mean sense is proved for fixed forcing terms with a fixed Hessian subsample size, while superlinear convergence in the mean square sense is obtained for the forcing terms that approach zero and increasing Hessian subsample sizes.

Having in mind the binary classification problem and the fact that the number of training points is enlarged over time in many applications, we also consider the case of unbounded $ N, $ i.e., the case where the objective function is defined as the mathematical expectation.  For this problem we  obtain convergence in the mean square sense as well.

Finally, the strong convexity assumption is relaxed similarly to the problem considered in  \cite{ mit1,mit2}.   A bound on the Hessian sample size, derived in  \cite{xu2} is used to obtain Hessians approximations positive definite with some high probability and CG is adapted to deal with possibly singular problems. The convergence in the  the mean square  is obtained for this problem as well.

From the performed theoretical analysis we derive adaptive rules for selecting  both the forcing term sequence and the Hessian sample size.
 Particularly important feature of the proposed method   is that the Hessian sample size is related to the current forcing term and approximated gradient norm, both quantities actually computable. Moreover, when $q$-superlinear convergence is sought, the Hessian sample size is adaptively chosen along  the iterations,  a low accuracy and smaller Hessian sample size are generally used in the early stage of the method while the accuracy and the Hessian sample size increase in the last stage of the convergence.  Finally, we note that the Hessian sample size is also allowed to  decrease if too high accuracy is used at the previous  iterations.
Numerical results on binary classification problems  give numerical evidence of the proposed adaptive choices effectiveness.

The paper is organized as follows. In Section 2 the method is introduced and the global and local analysis is performed using the standard deterministic reasoning for the case of strongly convex functions $ f_i$.    Convergence in the mean square is proved in Section 3, considering all three cases - finite number of strongly convex functions $ f_i, $ an infinite number of strongly convex functions $ f_i$'s and the last case with relaxed convexity assumptions. Some numerical results are presented in Section 4.

\subsection{Related work}

Our analysis is strictly related to that developed in  \cite{bbn_2017,bbn, mit1,mit2}, where convergence of Inexact subsampled Newton methods is investigated both in probability   \cite{ mit1,mit2} and expectation \cite{noc2, bbn, bbn_2017}.   We differ from these papers as  we focus  on the choice of the forcing terms,  on the nonmonotone  line search strategy and on adaptive choices of Hessian sample size.

Results presented in \cite{bbn,mit2} focus on the local behaviour of the method and give bounds for the accuracy required in the last stage of the procedure. The issue of developing an automatic transition between the initial stage   of the procedure where
a low accuracy in the Hessian approximation is enough, and the last stage stage where more accurate Hessian approximations are needed, is not investigated. However, in \cite{mit2} the analysis is carried out  under weaker assumptions than those we used here as the function $f_{\cal N}$ is supposed to be strongly convex only in a neighborhood of the sought  minimizer without any assumptions on convexity of functions $f_i$.   A set of conditions on the gradient and the Hessian sample sizes that ensures  local R-superlinear convergence in the expectation under the  assumption on the variance of the error norms (Bounded moment of Iterates) is given in  \cite{bbn}. The Hessian sample size is assumed to be increased at each iteration starting from a large enough initial sample size.

 In  \cite{noc2} an adaptive rule for  choosing the gradient sample size is proposed along with an automatic criterion for the forcing term related to a variance estimation of the Hessian accuracy.
The Hessian sample size is a fixed fraction of of the used gradient sample size.

Finally, we would like to mention that in  \cite{ bbn_2017}  the authors perform a local complexity analysis of subsampled Inexact Newton methods and  also show that methods that incorporate stochastic second-order information can be far more efficient on badly-scaled or ill-conditioned problems than first-order methods.

\section{Inexact subsampled Newton method}

 We first introduce the notation and give some preliminary results.
 Throughout the paper
$ {\cal N}_k\subset \{1,2,\ldots,N\} $ denotes the sample used to approximate the objective function and its gradient,
$ N_k $  denotes its cardinality and the subsampled function and gradient are defined as
$$  f_{{\cal N}_k}(x) = \frac{1}{N_k} \sum_{i\in {\cal N}_k } f_i(x), \;\;\;\; \nabla f_{{\cal N}_k} = \frac{1}{N_k} \sum_{i\in {\cal N}_k} \nabla f_i(x). $$  Moreover,  $ {\cal D}_k \subset \{1,2,\ldots,N\}$ is the sample used for Hessian approximation with cardinality $ D_k, $ and the subsampled Hessian is given by
\be \label{subhessian}   \nabla^2 f_{{\cal D}_k} (x)= \frac{1}{D_k} \sum_{i \in {\cal D}_k} \nabla^2 f_i(x). \ee

Here, we will consider subsampled  Inexact Newton methods, that is  iterative processes where at  iteration $k$, given the current iterate $x^k$, the step $s^k$ used to update the iterate satisfies
\be \label{inexact1}
 \nabla^2 f_{{\cal D}_k}(x^k) s^k = - \nabla f_{{\cal N}_k} (x^k) + r^k, \; \|r^k\| \leq \eta_k \|\nabla f_{{\cal N}_k}(x^k)\|.
\ee
The term $\eta_k$  belongs to $(0,1)$ and it is called forcing term  \cite{DES,nash}. Here and in the rest of paper  $\|\cdot\|$ denotes the 2-norm.

Through the paper we will restrict our attention to convex functions, more precisely we will first consider strongly convex functions and then, in Subsection 3.2 relax the strong convexity to convexity. Let us state this formally.

\noindent
{\bf Assumption A1} The functions
 $ f_i, \; i=1,\ldots,N $ are  twice continuously differentiable and strongly convex, i.e.,
\be \label{convexity}
\lambda_1 I \preceq \nabla^2f_i(x) \preceq \lambda_n I,\;\;\; \forall x \in {\mathbb R}^n\;\;\;\;i=1,\ldots,N
\ee
with $ \lambda_1 > 0. $

 Assumption A1 implies a couple of inequalities that will be used further on. First of all,
 for  all $ x \in \mathbb{R}^n $ we have
\be \label{Lemma2}  \lambda_{1} \|x-x^*\| \leq \|\nabla f_{\cal N}(x) \| \leq \lambda_{n} \|x-x^*\|,
\ee
 where $ x^* $ is the unique minimizer of the function $f_{{\cal N}}$.
Furthermore, according to \cite[Theorem 2.10]{nesterov},  for every $x \in \mathbb{R}^{n}$
\label{app1}
\be \label{conv1}  \frac{\lambda_1}{2 } \|x-x^{*}\|^{2}\leq f_{{\cal N}}(x)-f_{{\cal N}}(x^{*})\leq \frac{1}{\lambda_1} \|\nabla f_{{\cal N}}(x)\|^{2}. \ee

Since $  \nabla^2 f_{{\cal D}}(x) $ is positive definite we  choose CG as the linear solver for computing $s^k$ in  (\ref{inexact1}). Thus,  we assume that CG initialized with zero vector is employed at each Inexact Newton iteration. We will make use of the following technical Lemma that it is proved in \cite{fg}.

\begin{lem} \cite{fg} \label{cglemma1} Let $ Ax = b, $ where $ A $ is a symmetric and positive definite matrix. Furthermore, let us assume that CG is applied to this system and it is terminated  at the $i$th iteration. Then if CG is initialized with the null vector the approximate solution $ x_i $ satisfies
$ x_i^T A x_i = x_i^T b. $
\end{lem}

Now, the above Lemma and Assumption A1 clearly imply the following result.
\begin{lem} \label{Lemma3a}
Assume that $ s^k $ satisfying (\ref{inexact1}) is obtained through CG method initialized with the null vector applied to the linear system
$$
 \nabla^2 f_{{\cal D}_k}(x^k) s= - \nabla f_{{\cal N}_k} (x^k).
 $$
Then, $ \|s^k\| \leq \lambda_1^{-1} \|\nabla f_{{\cal N}_k}(x^k)\|. $
\end{lem}

%
%

\subsection{Global convergence}

In this section we analyze the behavior of the  subsampled IN method and CG as the inner solver wrapped into the nonmonotone line-search  strategy given in \cite{LF}. The choice of this line search rule is motivated by the fact that in the first stage of the procedure  we work with subsampled functions and gradients and enforcing a strict  Armijo type decrease might yield unnecessary small steps without real decrease in the original objective function.  Thus, we allow additional freedom in the step size selection introducing in  the Armijo condition a positive term usually denoted as error term, \cite{LF}. The line search procedure is governed by an error sequence $\{\nu_k\}$ with the following properties
\be \label{errorsequence}  \nu_k > 0, \quad \sum_k^\infty \nu_k < \infty. \ee
Iteration $k$ of the  above procedure, denoted as   Algorithm GIN, is sketched in Algorithm \ref{GIN}.

 \algo{GIN}{$k$-th iteration of Algorithm GIN}{
\noindent
Given  $ x^k\in {\mathbb R}^n \; c \in (0,1)$, $\bar \eta \in (0,1)$  and $ \{\nu_k\} $ such that (\ref{errorsequence}) holds.
\vskip 2 pt
\begin{itemize}
\item[Step 1.] Choose $ {\cal N}_k $, $ {\cal D}_k, $  $\eta_k\in (0,\bar \eta)$.
\item[Step 2.] Apply  CG method initialized by the null vector to  $\nabla^2 f_{{\cal D}_k}(x^k) s^k=- \nabla f_{{\cal N}_k}(x^k)$ and compute
$s^k$ satisfying (\ref{inexact1}).
\item[Step 3.] Find the smallest nonnegative integer $j$ such that for  $ t_k= 2^{-j} $ there holds
\be \label{lsnm}
f_{{\cal N}_k}(x^k + t_k s^k) \leq f_{{\cal N}_k}(x^k) + c t_k (s^k)^T \nabla f_{{\cal N}_k}(x^k) + \nu_k
\ee
and set $ x^{k+1} = x^{k} + t_k s^k$,  $ k = k+1$.
\end{itemize}
}\label{GIN}

The Algorithm GIN is stated with an arbitrary scheduling sequences $ \{N_k\} $ and $ \{D_k\} $.   Notice that here we define the line search rule with $ f_{{\cal N}_k} $ and thus inexact function values, as well as approximated  gradient values, are allowed in Algorithm GIN. Naturally, one aims at reaching the full sample at some iteration but expects to save computational effort while working with smaller samples. We are going to analyse the complexity and global convergence of this Algorithm. First of all, we will give iteration complexity result for an arbitrary schedule i.e. we will consider the decrease in (possibly inexact) gradient. Then, we will prove R-linear convergence of the iterative sequence
and show that the classical complexity result of  $ {\cal O}(\log(\epsilon^{-1})), $ \cite{grapiglia} can be proved for a schedule that eventually ends up with the full sample in a finite number of iterations.

As the search direction $ s^k $ is generated by the CG method and $ f_{{\cal N}_k} $ is strictly convex, we know by Lemma \ref{cglemma1} that $ s^k $ is  descent direction  for $ f_{{\cal N}_k}$ at $x^k$ and thus, under Assumption A1,  the step size  $ t_k $ is strictly positive even for the standard Armijo rule in Step 3. The lower bound for $ t_k $ is obtained in Lemma \ref{lowertk}, similarly to \cite{nas1}. 


\begin{lem} \label{cglemma2}
Let  $ s^k $ be the   step generated in Step 2 of Algorithm GIN.   Then
$$ - (\nabla f_{{\cal N}_k}(x^k))^T s^k \geq  \frac{\lambda_1}{\lambda_n^2} (1- \eta_k)^2 \|\nabla f_{{\cal N}_k} (x^k)\|^2. $$
\end{lem}

{\em Proof.} The inexact condition implies
$$
\|\nabla f_{{\cal N}_k}(x^k) - r^k\| \geq \|\nabla f_{{\cal N}_k}(x^k)\| - \|r^k\| \geq (1-\eta_k) \|\nabla f_{{\cal N}_k} (x^k)\|.
$$
On the other hand, we have
$$
\|\nabla f_{{\cal N}_k}(x^k) - r^k\| = \|\nabla^2 f_{{\cal D}_k}(x^k) s^k\| \leq \|\nabla^2 f_{{\cal D}_k}(x^k)\| \| s^k\| \leq \lambda_n \|s^k\|
$$ due to  (\ref{convexity}). Therefore
\be \label{in1}
\|s^k\| \geq \frac{\|\nabla f_{{\cal N}_k}(x^k) - r^k\|}{\lambda_n} \geq \frac{1-\eta_k}{\lambda_n} \|\nabla f_{{\cal N}_k}(x^k)\|.
\ee
Then, Lemma \ref{cglemma1} and (\ref{convexity}) yield
\begin{eqnarray}
- (\nabla f_{{\cal N}_k}(x^k))^T s^k & = & (s^k)^T \nabla^2 f_{{\cal D}_k}(x^k) s^k \geq \lambda_1 \|s^k\|^2   \label{bound_prod_a}\\
& \geq & \frac{\lambda_1}{\lambda_n^2} (1- \eta_k)^2 \|\nabla f_{{\cal N}_k}(x^k)\|^2. \label{bound_prod}
\end{eqnarray}
$\Box$

 \begin{lem}\label{lowertk}
Suppose that the assumption A1 holds and let $ s^k $ be  generated in Step 2 of algorithm GIN. Then (\ref{lsnm}) holds for  $t_k \geq \bar{t}= (1-c)\lambda_1/\lambda_n. $
\end{lem}

{\em Proof.}
 Let $ k $ be an arbitrary iteration. 
 If $t_k=1$ satisfies \eqref{lsnm}, that is  in Step 3 we have $j=0$, then $t_k$ is  greater than $\bar t$, as $\bar t\in (0,1)$. 
 So let us consider the case  $t_{k}<1$. Then there exists $t'_{k} = 2 t_{k} $ such that
\begin{eqnarray*}
f_{{\cal N}_k}(x^{k}+t'_{k} s^{k})&>& f_{{\cal N}_k}(x^{k})+c t'_{k} (s^{k})^{T} \nabla f_{{\cal N}_k} (x^{k})+\nu_k\\& \geq &  f_{{\cal N}_k}(x^{k})+ c t'_{k} (s^{k})^{T} \nabla f_{{\cal N}_k} (x^{k}) .
\end{eqnarray*}
On the other hand, Assumption A1 implies, using the standard arguments for functions with bounded Hessians,
\begin{eqnarray*}
f_{{\cal N}_k}(x^{k}+t'_{k} s^{k})&=&
f_{{\cal N}_k}(x^{k}) +  \int_{0}^{1} (\nabla f_{{\cal N}_k}(x^{k}+y t'_{k} s^{k}))^{T} ( t'_{k}s^{k}) dy\\
&\leq & \frac{\lambda_n}{2} (t'_{k})^{2}\| s^{k}\|^{2}+ f_{{\cal N}_k}(x^{k})+t'_{k}(\nabla f_{{\cal N}_k}(x^{k}))^{T}s^{k}.\\
\end{eqnarray*}
Combining the previous two inequalities we obtain
$$c t'_{k} (s^{k})^{T} \nabla f_{{\cal N}_k}(x^{k}) \le \frac{\lambda_n}{2} (t'_{k})^{2}\| s^{k}\|^{2}+ t'_{k}(\nabla f_{{\cal N}_k}(x^{k}))^{T}s^{k}.$$
Dividing by $t'_{k}$ and using  $t_{k}=t'_{k} /2$, by rearranging the previous inequality we get
\be \label{rlin6} t_k \geq \frac{-(1-c)(\nabla f_{{\cal N}_k}(x^{k}))^{T}s^{k}  }{\lambda_n \| s^{k}\|^{2}}.\ee
Now, the result follows from  \eqref{bound_prod_a} and the fact
$ \min\{1,(1-c)\lambda_1/\lambda_n\} = (1-c)\lambda_1/\lambda_n. $
$ \Box$

To prove the main results we need the following Lemma from \cite{nas1}.

\begin{lem} \label{mnkj} Assume that $ \zeta_k \to 0 $ R-linearly.  Then, for every $\rho \in (0,1)$,
$$a_{k}=\sum_{j=1}^{k}\rho^{j-1} \zeta_{k-j}$$
converges to zero R-linearly.
\end{lem}


Let us denote with $\xi^g_k $ and  $ \xi^f_k $  the inaccuracy in function and gradient, i.e.   $\xi^g_k $ and  $ \xi^f_k $ are such that:
\be\label{errfg}
 \max_{x \in \{x^k, x^{k+1}\}} |f_{{\cal N}_k}(x) -f_{{\cal N}}(x)| \le \xi^f_k,  \quad
 |\|\nabla f_{{\cal N}_k}(x^k)\|^2 -\|\nabla f_{{\cal N}}(x^k)\|^2 |\le \xi^g_k.
\ee
Following Grapiglia and Sachs \cite{grapiglia}, we now prove that despite inaccuracy in function and gradient, Algorithm GIN meets the complexity results of nonmonotone line search methods with exact function and gradients.

 \begin{teo} \label{complexity}
 Suppose that Assumption A1 holds and
\be \label{errorfunction}   \sum_{k=0}^\infty \xi_k^f < \infty. \ee
Then, for a given $\epsilon \in (0,1)$ Algorithm  GIN takes at most
$$ \bar k=\lceil \frac{f_{\cal N}(x^0)- f_{\cal N}(x^*)+ \sum_{k=0}^\infty (2 \xi_k^f +\nu_k)}{\kappa_c}\epsilon^{-2}\rceil,
$$
iterations to ensure $\|\nabla f_{{\cal N}_{\bar k}}(x^{\bar k}) \|\le \epsilon$, where
\be\label{kc}
\kappa_c =
c (1-c)  \frac{\lambda_1^2}{{\lambda_n}^3} (1 -\bar{\eta})^2.
\ee
\end{teo}
{\em Proof.}
Note that,  by Lemma  \ref{cglemma2} there follows
\begin{equation}\label{compl}
 (\nabla f_{{\cal N}_k}(x^k))^T s^k \leq  -\frac{\lambda_1}{\lambda_n^2} (1- \bar \eta)^2 \|\nabla f_{{\cal N}_k}(x^k)\|^2.
\end {equation}
Then, we can proceed as in the proof of Theorem 1 in \cite{grapiglia}.
Let $\bar k$ be the first iteration such that $\|\nabla f_{{\cal N}_{\bar k}}(x^{\bar k}) \|\le \epsilon$.
By \eqref{lsnm}  we obtain
$$
\nu_k+f_{{\cal N}_k}(x^k)-f_{{\cal N}_k}(x^{k+1})\geq -c t_k (\nabla f_{{\cal N}_k}(x^k))^T s^k, \; k=0,1,\ldots,\bar{k}-1.
$$
 Moreover, by \eqref{compl} and  Lemma  \ref{lowertk}, there follows
$$-c t_k (\nabla f_{{\cal N}_k}(x^k))^T s^k  \geq c (1-c)  \frac{\lambda_1^2}{{\lambda_n}^3} (1 -\bar{\eta})^2 \|\nabla f_{{\cal N}_k}(x^k)\|^2.
$$
Then, for $ k=0,1,\ldots, \bar{k}-1, $ there holds
$$ \nu_k+f_{{\cal N}_k}(x^k)-f_{{\cal N}_k}(x^{k+1})\geq  \kappa_c \epsilon^2 $$
with $ \kappa_c$ given in \eqref{kc}.

Therefore, by \eqref{errfg}
$$
\nu_k+2 \xi_k^f+ f_{{\cal N}}(x^k)-f_{{\cal N}}(x^{k+1})\ge \kappa_c   \epsilon^2 \quad k=0,\ldots,\bar k-1.
$$
Summing up for $k=0,\ldots,\bar k-1$ we get
$$
\sum_{k=0}^\infty (2 \xi_k^f +\nu_k)+f_{\cal N}(x^0)-f_{{\cal N}}(x^{\bar k})\ge \bar k \kappa_c \epsilon^2
$$
and the thesis follows. $\Box$

 Note that in the previous theorem if $N_{\bar k} <N$ then
\be\label{conv_grad}
\|\nabla f_{\cal N}(x^k) \|\le \epsilon _{\bar k}+\xi^g_k.
\ee


 In what follows we prove $R$-linear convergence and consequently complexity estimates which are far
better than $O(\epsilon^{-2})$ at the price of a strengthened assumption on the accuracy in functions and gradients.

\begin{teo} \label{thglobal}
Assume that A1 holds and  let $ \{x^k\} $ be generated by Algorithm GIN. If the error sequences $ \{\nu_k\} $,
$ \{\xi_k^f\} $ and $ \{\xi_k^g\} $
 converge to zero R-linearly then  $ \{x^k\} $ converges R -linearly to the solution of (\ref{problem1}).
\end{teo}

{\em Proof.} 
Inequalities (\ref{lsnm}), (\ref{errfg}) and (\ref{conv1}) and Lemma \ref{cglemma2} imply
\begin{eqnarray*}
\fnk(x^{k+1}) - f_{\cal N}(x^*) & \leq & \fnk(x^k) - f_{\cal N}(x^*) + c t_k \nabla \fnk(x^k)^T s^k + \nu_k\\
& \leq &  \fnk(x^k) - f_{\cal N}(x^*) - c t_k \frac{\lambda_1}{\lambda_n^2} (1- \eta_k)^2 \|\nabla \fnk (x^k)\|^2 + \nu_k\\
& \leq & \fnk(x^k)-f_{\cal N}(x^*) \\
& &- c t_k \frac{\lambda_1}{ \lambda_n^2} (1-\eta_k)^2(\lambda_1 (f_{\cal N}(x^k)-f_{\cal N}(x^*))-\xi_k^g)+\nu_k
\end{eqnarray*}
Then, using  Lemma \ref{lowertk} and (\ref{errfg}) again, we obtain
\begin{eqnarray} \label{novo}
f_{\cal N}(x^{k+1}) - f_{\cal N}(x^*)
& \leq &\rho (f_{\cal N}(x^k)-f_{\cal N}(x^*)) + \bar \xi_k
\end{eqnarray}
where $\rho=1- c \bar t \frac{\lambda_1^2}{ \lambda_n^2} (1-\bar\eta)^2 \in (0,1)$ and $\bar \xi_k=\nu_k+2\xi_k^f+\xi_k^g$. Furthermore, we obtain
\be \label{Rlinear} f_{\cal N}(x^{k+1}) - f_{\cal N}(x^*)\leq \rho^{k+1} (f_{\cal N}(x^0) - f_{\cal N}(x^*)) + \rho \sum_{j=1}^{k} \rho^{j-1} \bar \xi_{k-j} + \bar \xi_{k}.\ee
Thus, Lemma \ref{mnkj} yields the statement.
$\Box$

Notice that R-linear convergence result obtained in Theorem \ref{thglobal} also holds for $\nu_k=0$, that is, for the Armijo line search.
The theorem above also allows us to prove the complexity result below.

\begin{teo} \label{compl2}
Assume that A1 holds.  If the error sequences $ \{\nu_k\} $,
$ \{\xi_k^f\} $ and $ \{\xi_k^g\} $
 converge to zero R-linearly then  for any $ \epsilon \in (0, e^{-1}) $,  there exist $\hat{\rho} \in (0,1)$ and $Q > 0$ such that Algorithm GIN  takes at most
 $$
\bar k =  \lceil \frac{\log((f_{\cal N}(x^0)- f_{\cal N}(x^*)+Q))}{|\log(\hat{\rho})|}\rceil \log(\epsilon^{-1})
$$
iterations to ensure $ f_{{\cal N}}(x^{\bar k}) -   f_{{\cal N}}(x^*) < \epsilon $.
\end{teo}

{\em Proof.}  Assumptions of Theorem \ref{thglobal} are satisfied and
\be \label{Rlinear1}
 f_{\cal N}(x^{k}) - f_{\cal N}(x^*)  \leq   \rho^k (f_{\cal N}(x^{0}) - f_{\cal N}(x^*)) + \sum_{j=0}^{k-1}  \rho^{j} \bar \xi_{k-j}. 
 \ee
 Given that the sequence $\{\sum_{j=0}^{k} \rho^{j} \bar \xi_{k-j}\} $ converges to zero R-linearly by Lemma \ref{mnkj}, there exist $ \bar{\rho} \in (0,1) $ and $ Q > 0 $ such that
$$ \sum_{j=0}^{k-1}  \rho^{j} \bar \xi_{k-j} \leq Q \bar{\rho}^{k}. $$
 Therefore, for  $\hat{\rho}=\max \{\rho, \bar{\rho}\},$ we have
from \eqref{Rlinear1}
$$
 f_{\cal N}(x^{k}) - f_{\cal N}(x^*)  \leq  \hat{\rho}^k  (f_{\cal N}(x^0) - f_{\cal N}(x^*)+Q)
$$
and the statement follows as in Theorem 6 of \cite{grapiglia}. $ \Box$

A couple of comment is due here. First of all, Theorems \ref{complexity}-\ref{compl2} deal with the possibility of infinite sequence of errors in the objective function and the gradient. So, these statements provide a framework for considering unbounded $N $ as well. However, our focus here is on the problems with fixed $ N. $ Thus,  any kind of scheduling for $ N_k $ that ensures the existence of $ k_0 $ such that $ N_k = N $ for $ k \geq \tilde{k}$ actually implies that the sequences $ \{\xi_k^f\} $ and $ \{\xi_k^g\} $ are  finite, with $ \tilde{k} $ nonzero elements at most. Therefore, the statements of the above  theorems apply to any kind of scheduling that ensures reaching the full sample for the objective function. Theorem \ref{compl2} proves the complexity bound of $ \log(\epsilon^{-1}), $ although we work with cheaper objective function and the gradient whenever $ N_k < N. $ Thus the theorem above provides theoretical justification for working with smaller samples.
{Moreover, the analysis carried out so far does not involve either the forcing term or the accuracy in Hessian approximation. Then, the results we provided hold even if $D_k=1$ and only one CG-iteration is performed at each iteration of Algorithm GIN.

\subsection{Local convergence}
 In this section we assume that the scheduling of the sample sizes $ {\cal N}_k $ is given and that eventually we reach $ f_{\cal N} $ and $ \nabla f_{\cal N} $ at some iteration $ \tilde k $ and continue with $ {\cal N}_k = {\cal N} $ for $ k \geq \tilde k.$ Then, we may restrict all asymptotic theoretical considerations to  the method defined with   ${\cal N}_k={\cal N}$, that is the step $s^k$ satisfies
\be \label{inexact4}
 \nabla^2 f_{{\cal D}_k}(x^k) s^k = - \nabla f_{\cal N} (x^k) + r^k, \; \|r^k\| \leq \eta_k \|\nabla f_{\cal N}(x^k)\|.
\ee
Following the analysis in \cite{EW1}, in this subsection  we focus on  the local $q$-linear and $q$-superlinear convergence of  method GIN. Then, we assume that the generated sequence converges to the solution of \eqref{problem1}, relying on the global convergence results proved in the previous subsection.
The analysis presented here differs from that in \cite{EW1} as we have to take into account the approximation in the Hessian and enlightens that the accuracy in the Hessian's approximation must be related to the adopted forcing term.
 This theoretical study enables us to devise an adaptive and computable rule for  the Hessian sample size yielding  $q$-linear and $q$-superlinear convergence. In fact, this is what we should obtain in order to justify the use of a second order method.
The analysis is carried out under  the following assumption on the subsampled Hessian matrices.

\vskip 5 pt
\noindent
{\bf Assumption A2.} There exists a constant $L$ such that
for any $ {\cal D} \subseteq \{1,2,\ldots,N\}  $ we have
\be \label{A4} \|\nabla^2 f_{\cal D}(x) - \nabla^2 f_{\cal D}(x^*)\| \leq L \|x - x^*\|,  \ee
for $x$ sufficiently close to $x^*$.

The above Assumption holds if  there exists a neighborhood  of $x^*$  where each Hessian $\nabla^2 f_i$
is Lipschitz continuous.
Note also that Assumption A1 implies
that all function $ \nabla f_i $ are Lipschitz continuous with the constant $ \lambda_n  $.

Moreover,  we assume that the error in the Hessian approximation   is determined only by the subsample size, independently of the sample taken and
denote by  $h({D,x})$ the norm of the error in the Hessian approximation for a given subsample size $D$  at point $ x$, i.e.
$$ h({D,x}) := \|\nabla^2 f_{{\cal N}}(x) - \nabla^2 f_{\cal D}(x)\|.$$
%
We make the following assumption on the subsample size $D_k$ used at each iteration $k$.

\vskip 5 pt
\noindent
{\bf Assumption A3.} At each iteration $ k$  of Algorithm GIN the subsample size    $ { D}_k $ is chosen such that
\be \label{Dk}
h({ D}_k, x^k) \leq C \eta_k,
\ee
 with $0 < C <(1/{\bar{\eta}} - 1)\lambda_1/2$.


 In the subsequent  analysis for any $\delta>0$  the ball with center $x^*$ and radius $\delta$ will be denoted by $ N_{\delta}(x^*)$, i.e.
$ N_{\delta}(x^*) = \{x \in \mathbb{R}^n: \|x - x^*\| \leq \delta\}.$ Moreover, we let
 $ \delta^*$ sufficiently small such that   (\ref{A4})   holds for any $x$ in $ N_{\delta^*}(x^*)$  and
{\bf$\delta^*< 2 \min \{\lambda_1, 1/ \lambda_n\}/L.$}

The next Theorem states that the full step is taken in Algorithm GIN eventually.


\begin{teo} \label{thRrate}
Assume that the sequence  $\{x^k\}$ generated by Algorithm GIN converges to $x^*$. Let Assumptions A1 and A3  hold  and $c\in (0,1/4)$ in (\ref{lsnm}).  Then,  there exists $ k_0 $ such that for all $ k \geq k_0 $ the full step $ s^k $  is accepted in Step 3 of Algorithm GIN.
\end{teo}

{\em Proof.}  Take $\bar  \varepsilon\in (0,\bar \eta\lambda_1/2) $ and $ \bar \delta >0$ such that
\be \label{ineq1}
\| \nabla^2 f_{\cal N} (x^k+ \xi s^k) - \nabla^2 f_{\cal N} (x^k) \| \leq \bar \varepsilon, \mbox{ for all } \|s^k\| \leq \bar \delta, \xi \in (0,1).
\ee
According to  Lemma \ref{Lemma3a},  $ \lim_{k \to \infty} \|s^k\| = 0$. Then, there exists $k_0$ such that   \eqref{ineq1} holds for $ k \geq k_0 $.
 Moreover, the Taylor expansion yields
\begin{eqnarray}
f_{\cal N}(x^k + s^k) & = & f_{\cal N}(x^k) + (\nabla f_{\cal N}(x^k))^T s^k + \frac{1}{2} (s^k)^T \nabla^2 f_{\cal N}(\theta^k) s^k \nonumber \\
& = &  f_{\cal N}(x^k) + (\nabla f_{\cal N}(x^k))^T s^k + \frac{1}{2} (s^k)^T \nabla^2 f_{{\cal D}_k}(x^k) s^k \nonumber \\
& &+  \frac{1}{2} (s^k)^T (\nabla^2 f_{\cal N}(\theta^k) -  \nabla^2 f_{\cal N}(x^k)) s^k  + \nonumber \\
& &+  \frac{1}{2} (s^k)^T (\nabla^2 f_{\cal N}(x^k) -  \nabla^2 f_{{\cal D}_k}(x^k)) s^k ,\label{ineq2}
\end{eqnarray}
with $ \theta^k = x^k + \xi s^k, \; \xi \in (0,1). $
 From Assumption A3 we obtain
\be \label{ineq3}
(s^k)^T (\nabla^2 f_{\cal N}(x^k) -  \nabla^2 f_{{\cal D}_k}(x^k)) s^k \leq h(D_k,x^k) \|s^k\|^2 \leq C \eta_k \|s^k\|^2 \leq C \bar \eta \|s^k\|^2.
\ee
Recall that Lemma \ref{cglemma1} implies
\be \label{ineq4}
(s^k)^T \nabla^2 f_{{\cal D}_k}(x^k) s^k = - \nabla f_{\cal N}(x^k)^T s^k.
\ee
Putting  (\ref{ineq1}), (\ref{ineq3}) and (\ref{ineq4}) into (\ref{ineq2}) we get
\be \label{ineq5}
f_{\cal N}(x^k + s^k)  \leq  f_{\cal N}(x^k) + \frac{1}{2} (\nabla f_{\cal N}(x^k))^T s^k + \frac{1}{2} (\varepsilon + C \bar{\eta}) \|s^k\|^2.
\ee
From \eqref{bound_prod_a}, we conclude that
$$ f_{\cal N}(x^k + s^k)  \leq  f_{\cal N}(x^k) + \frac{1}{2}(1-\frac{\varepsilon + C \bar{\eta}}{\lambda_1}) ( \nabla f_{\cal N}(x^k))^T s^k.$$
Note that from   $ C< (1/\bar{\eta}-1)\lambda_1/2 $ there follows $\frac{C\bar\eta}{\lambda_1}<(1-\bar \eta)/2$. Therefore, the choice of $\varepsilon$ yields
$ \varepsilon/\lambda_1 + C\bar\eta/\lambda_1 < 1/2$. Then,
$$ f_{\cal N}(x^k + s^k)  \leq  f_{\cal N}(x^k) + \frac{1}{4} \nabla f_{\cal N}^T(x^k) s^k$$
and condition (\ref{lsnm})
is satisfied with  $t_k=1$ for any $k>k_0$ as $c\in (0,1/4)$.   $ \Box $



The following  Lemma, whose proofs can be found in the Appendix, is needed in the subsequent convergence analysis.

\begin{lem} \label{Lemma4}
 Let  Assumptions A1-A2 hold.
If $ x^k \in {\cal N}_{\delta^*}(x^*), $ $ s^k \in \mathbb{R}^n $ and $ \eta \in (0,1) $ are such that $  \|\nabla f_{\cal N}(x^k) + \nabla^2 f_{{\cal D}_k}(x^k) s^k \| \leq \eta \|\nabla f_{\cal N}(x^k)\| $ and $ x^k + s^k \in {\cal N}_{\delta^*}(x^*) $ then $$ \|\nabla f_{\cal N}(x^k+s^k)\| \leq (\eta + B(x^k) ) \|\nabla f_{\cal N}(x^k)\|, $$ with $ B(x^k)  =\frac{1}{\lambda_1}  (\frac{1}{2\lambda_1} L \|\nabla f_{\cal N}(x^k)\| + h(D,x^k)). $
\end{lem}

\begin{lem} \label{Lemma5}
Let  Assumption A1 holds and $\delta\in(0,\delta^*/(1+\lambda_1^{-1} \lambda_n))$. If $ x^k \in {\cal N}_{\delta}(x^*) $ and $ \|\nabla f_{\cal N}(x^k) + \nabla^2 f_{{\cal D}_k}(x^k) s^k\| \leq \eta \|\nabla f_{\cal N}(x^k)\| $ then $ x^k+s^k \in N_{\delta^*}(x^*). $
\end{lem}

The convergence of $ \{x^k\} $ together with  (\ref{Lemma2}),   implies the following result.

\begin{teo} \label{Cor2}
Assume that the sequence  $\{x^k\}$ generated by Algorithm GIN converges to $x^*$. Let Assumptions A1-A3 hold and $c\in (0,1/4)$ in (\ref{lsnm}).  If  $\eta_k = \bar \eta$  at each iteration of Algorithm GIN,  then the sequence $ \{x^k\} $ converges to $ x^* $ q-linearly for $\bar  \eta $ small enough. Moreover, if   $ \lim_{k \to \infty} \eta_k = 0, $
the convergence is $q$-superlinear.
\end{teo}
{\em Proof.}
Note that, by the choice of $C$ in Assumption A3, $ \bar \eta(1+\lambda_1^{-1} C)<(1+\bar \eta)/2$.  Let
$\delta\in (0,\delta^*/(1+\lambda_1^{-1} \lambda_n))$.
Take $\bar  \varepsilon \in ( 0,  \delta \lambda_1]$ sufficiently small such that $ \bar \eta(1+\lambda_1^{-1} C) + \lambda_1^{-2} L \bar \varepsilon/2 < \tau <1 $ for some $ \tau \in (0,1)$.
Let $k_0$ be  defined as in Theorem \ref{thRrate} and $\bar k\ge k_0$, such that
$ x^{\bar k}\in   {\cal N}_{\delta}(x^*)$   sufficiently near to $x^*$ to guarantee $ \|\nabla f_{\cal N}(x^{\bar k})\| \leq \bar \varepsilon. $

 Lemma \ref{Lemma5} yields $ x^{\bar k+1} \in {\cal N}_{\delta^*}(x^*) $ and  by Lemma \ref{Lemma4} we obtain
$$ \|\nabla f_{\cal N}(x^{\bar k+1})\| \leq \tau \|\nabla f_{\cal N}(x^{\bar k})\| \leq \|\nabla f_{\cal N}(x^{\bar k})\| \leq \bar\varepsilon. $$  Therefore, using (\ref{Lemma2}),
$$ \|x^{\bar k+1} - x ^*\| \leq \frac{1}{\lambda_1} \|\nabla f_{\cal N}(x^{\bar k+1})\| \leq \frac{1}{\lambda_1} \bar \varepsilon \leq \delta, $$ so $ x^{\bar k+1} \in {\cal N}_{\delta}(x^*)$ and  $ \|\nabla f_{\cal N}(x^{\bar k+1})\|\le \bar \varepsilon$.

As an inductive hypothesis suppose that for some $ k >{\bar k} $ we have $ x^k \in {\cal N}_{\delta}(x^*)$ and $ \|\nabla f_{\cal N}(x^k)\| \leq \bar \varepsilon $. Then $ x^{k+1} =x^k + s^k \in {\cal N}_{\delta^*}(x^*) $ by Lemma \ref{Lemma5}, and  Lemma \ref{Lemma4}   implies
\begin{eqnarray*}
  \|\nabla f_{\cal N}(x^{k+1})\| & \leq & [(1+\lambda_1^{-1} C) \eta_k + \lambda_1^{-2}  L \varepsilon/2] \|\nabla f_{\cal N}(x^k)\| \\
  & \leq & \tau \|\nabla f_{\cal N}(x^k)\| \leq \|\nabla f_{\cal N}(x^k)\| \leq \varepsilon.
\end{eqnarray*}
Again, (\ref{Lemma2})  yields $ \|x^{k+1} - x^*\| \leq \delta $
and $ x^{k+1} \in {\cal N}_{\delta}(x^*). $
Therefore,  proceeding by induction   we conclude that     $x^{k} \in {\cal N}_{\delta^*}(x^*)$ for any $k\ge \bar k$ and
 by  Lemma \ref{Lemma4},
\be \label{rate_ineq}
  \|\nabla f_{\cal N}(x^{k+1})\|  \leq  [(1+\lambda_1^{-1} C) \eta_k + \lambda_1^{-2}  L \|\nabla f_{\cal N}(x^k)\|/2] \|\nabla f_{\cal N}(x^k)\|,\quad k\ge \bar k.
\ee
Therefore, as $\|\nabla f_{\cal N}(x^k)\|\rightarrow 0$, using  (\ref{Lemma2}) we obtain  that
$ \{ x^k\}$ converges to $ x^* $ with  $q$-linear rate provided that
$$ \bar \eta<\frac{\lambda_1}{\lambda_n} \frac{1}{1+\lambda_1^{-1} C}<\frac{\lambda_1}{\lambda_n}.$$
Moreover, the  $q$-superlinear convergence follows if $ \lim_{k \to \infty} \eta_k = 0. $
$ \Box $

 The above results are in line with the classical convergence theory of Inexact Newton methods  \cite{DES} as the local linear convergence require $\eta_k\le \bar \eta <1$ and the upper bound on $\bar \eta$ depends on inverse of the conditioning of the Hessian.

 Let us now discuss one possible choice of $ \eta_k $ in order to obtain $q$-superlinear convergence of the procedure.  Following ideas in \cite{EW1} our choice of $\eta_k$ depends on the agreement between the function and the subsampled Newton model. If there is a good agreement between these two quantities, even if the quality of the approximation in the Hessian is lower than that in the function, it is reasonable to use a small $\eta$ in the subsequent iteration. Let us consider the following choice of $\eta_k$ in Algorithm \ref{GIN},
\be \label{eta2}  \eta_{k}=\min \{\bar{\eta},\frac{|f_{{\cal N}_k}(x^{k})-m_{k-1} (s^{k-1})|}{\|\nabla f_{{\cal N}_{k-1}}(x^{k-1})\|}\}, \; \bar{\eta}<1\ee
where
\be \label{genmod}  m_{k-1}(s) = f_{{\cal N}_{k-1}}(x^{k-1}) + \nabla f_{{\cal N}_{k-1}}(x^{k-1})^T s + \frac{1}{2} s^T \nabla^2 f_{{\cal D}_{k-1}}(x^{k-1}) s. \ee

In the following theorem we show that  the sequence  $ \{\eta_k\} $ generated by (\ref{eta2})   converges to zero and
ensures q-superlinear convergence of the sequence generated by Algorithm GIN if the scheduling ensures that the full sample is reached at some finite iteration.
\begin{teo} \label{Theorem3}
Let assumptions in Theorem \ref{Cor2}  hold and
 $ \eta_k $ given by (\ref{eta2}).  Then $ \{x^k\} $ converges to $ x^* $ superlinearly.
\end{teo}

{\em Proof.}
 Let the iteration index $k>\tilde k+1$, i.e. such that ${\cal N}_{k-1}={\cal N}_{k}={\cal N}$.
Using the Taylors expansion and (\ref{Dk}) we obtain
 \be \label{razlika} |f_{\cal N}(x^{k}) - m_{k-1}(s^{k-1})|\leq \frac{1}{2} \|s^{k-1}\|^2 \left(\frac{L}{2}\|s^{k-1}\|+C \eta_{k-1}\right).\ee
Now, by Lemma \ref{Lemma3a} there follows
\be \label{eta3}
\eta_k \leq \frac{1}{2} \lambda_1^{-2} \|\nabla f_{\cal N}(x^{k-1})\|^2\left[ \lambda_1^{-1} \frac{L}{2}  \|\nabla f_{\cal N}(x^{k-1})\| + C \bar{\eta}\right].
\ee
Then, as   $ \lim_{k \to \infty} \|\nabla f_{\cal N}(x^{k})\| = 0 $ we have   $ \lim_{k \to \infty} \eta_k = 0 $  and  the superlinear convergence follows by Theorem \ref{Cor2}. $ \Box $

 The above result can be proved also choosing $\eta_k$ as
\be \label{eta1}
\eta_k =\min \{\bar{\eta}, \frac{| f_{{\cal N}_k}(x^{k}) - m_{k-1}(s^{k-1})|}{\omega_k}\}, \; \bar{\eta}<1
\ee
with
 $$ \omega_k = {\cal O}(\|\nabla f_{{\cal N}_{k-1}} (x^{k-1})\|). $$
We also underline that Theorem  \ref{Theorem3} in case of  full Hessian, i.e. ${\cal N}_k={\cal D}_k
$,
shows  that the Truncated  Newton method for unconstrained optimization problems with the choice of forcing terms given by (\ref{eta1})  is superlinearly
convergent.     Such result also follows from the analysis presented in \cite[Theorem 3.10]{lee}.

\section{Mean square convergence}

Large part of the previous analysis strongly relays on the Hessian error bound $h(D,x)$. However, assessing this quantity is not an easy task in general. In this section, we provide some more specific guidance on choosing $D_{k}$. We assume that the subset ${\cal D}_{k}$ is chosen randomly and we relax controlling $h(D,x)$ in such way that we ask for a good enough Hessian approximation with some probability smaller than 1. Thus less conservative estimates are now feasible. More precisely, we use a bound similar to that derived in \cite{xu2} to carry out our analysis. Consequently, we obtain stochastic convergence results - convergence in a mean square sense (m.s.). Our main result is that the $q$-linear convergence in  m.s. can be achieved with a small enough but fixed forcing term $\eta$ and with  large enough but fixed Hessian sample size $D$. On the other hand, to achieve  $q$-superlinear convergence in  m.s. with $\eta_{k}$ defined by (\ref{eta1}),  the Hessian sample size is required to increase as $\eta_k$ goes to zero. This analysis paths the way to devise adaptive rules for selecting the Hessian sample size in order to use a
small Hessian sample size  in the early stage of the procedure, when $\eta_k$ is close to one and linear systems are solved only to a low accuracy, and automatically increase it when the solution is approached.

In this section we  assume that the subsample ${\cal D}$ is chosen randomly and uniformly - every $\nabla^2 f_i(x)$ has the same chance to be chosen. Let $ {\cal D} $ be any subset of $ {\cal N} $ such that $ |{\cal D}| = D. $ 
Then one can derive a bound on $D$ such that the following holds for a given $ \gamma > 0 $ and $ \alpha \in (0,1) $
\be \label{cilj} P(\|\nabla^2 f_{\cal D}(x)-\nabla^2 f_{\cal N}(x)\|\leq \gamma ) \geq 1-\alpha. \ee
The corresponding bound is stated in Lemma \ref{dkbounds} (see the proof in the Appendix), while  a similar bound is provided in Lemma 4 of \cite{xu2}. The result is obtained by using the Bernstein inequality, see \cite{tropp12} and \cite{bern} for further references.}

 \begin{lem} \label{dkbounds}
Assume that A1 holds and that the subsample ${\cal D}$ is chosen randomly and uniformly from ${\cal N}$. Let $ \gamma> 0 $ and $ \alpha \in (0,1) $ be given. Then (\ref{cilj}) holds at any point $x$ if the subsample size $ D$ satisfies
\be \label{minbound} D \geq \frac{2 (\ln (2n/ \alpha))(\lambda_n^2+ \lambda_n \gamma/3)}{\gamma^2}:=\tilde{l}.\ee
\end{lem}

 We use the above results and the analysis of the previous section, to design  a globally convergent inexact subsampled
Newton method with adaptive choice of the Hessian sample size. In particular we will choose $D_k$ such that the inequality
\begin{equation}\label{Dk_choice}
h(D_k,x^k)\le C \max \{\eta_{k}, \|\nabla f_{{\cal N}_k}(x^{k})\|\}
\end{equation}
holds
with probability $1-\alpha_k$, with $\alpha_k\in (0,1)$ and  $0 < C <(1/{\bar{\eta}} - 1)\lambda_1/2$.
This corresponds to  $\gamma= C \max \{\eta_{k}, \|\nabla f_{{\cal N}_k}(x^{k})\|\}$ in  (\ref{cilj}).

\subsection{Bounded sample - GIN-R method}

In this subsection we are interested in the case of finite $ N. $ Let us assume that the full sample is eventually reached for the objective function and the gradient, i.e., $ {\cal N}_k = {\cal N}, $ for $k$ sufficiently large, say $k\ge \bar k. $ The procedure we obtain  is  based on GIN method but with a specific choice of the Hessian subsample ${\cal D}_k$, namely its cardinality is set according to (\ref{Dk_choice})   and the sample is randomly chosen. We denote this procedure as GIN-R to emphasize the specific random choice of the Hessian subsample and for the sake of convenience we list  its generic iteration $k$ in Algorithm  \ref{GIN-R}, where we denote the first steps as Step 1.a-1.c to make clear that they correspond to specific choices in  Step 1 of Algorithm GIN.

We will analyze the convergence in the m.s.  considering two possibilities. If the forcing terms converge to zero, i.e., if $ \eta_k \to 0, $ then $ \gamma_k $ given by (\ref{62}) converges to zero. 
 The other case we  consider is $\eta_k=\bar{\eta} <1.$ In that case Algorithm GIN-R yields $\gamma_{k}$ bounded away from zero and we have $ {\cal D}_k \subset {\cal N} $ with $ D_k < N $ during the whole iterative process.
Theorem \ref{thglobal} implies the R-linear convergence for any Hessian subsampling under the condition  ${\cal N}_k = {\cal N}, k \geq \bar k.$    The key issue in the global convergence analysis is that we have  descent search direction with an arbitrary good or poor Hessian approximation, i.e., regardless of the subsample used in GIN algorithm. The line search globalization strategy makes the algorithm convergent. Therefore,  no matter how the subsample is chosen, the R-linear convergence still holds.
Thus we immediately have the following statement.

 \algo{GIN-R}{$k$-th iteration of Method GIN-R}{
\noindent
Given  $ x^k\in {\mathbb R}^n, \; \bar \eta \in (0,1), \;c \in (0,1), \;  C>0, \; \nu_k,\; \alpha_{k}\in (0,1).$
\vskip 5 pt
\noindent
Step 1.a Choose $ {\cal N}_k, \eta_k\in (0,\bar \eta). $
\vskip 5 pt
\noindent
Step 1.b Compute
\be \label{62} \gamma_{k}=C \max \{\eta_{k}, \|\nabla f_{{\cal N}_k}(x^{k})\|\}.\ee
\vskip 5 pt
\noindent
Step 1.c Set $\alpha=\alpha_k$ and $\gamma=\gamma_k$ and compute $ D$ such that (\ref{minbound}) holds. If
 $D \geq N_k$\\
\hspace*{33pt}  set ${\cal D}_k = {\cal N}_k$. Else, choose  the  sample $ {\cal D}_k $ randomly and uniformly from\\
  \hspace*{33pt}$ {\cal N}_k  $ such that $D_k \geq D. $
\vskip 5 pt
\noindent
Step 2 Apply  CG method initialized by the null vector to\\
\hspace*{33pt}  $\nabla^2 f_{{\cal D}_k}(x^k) s^k = -\nabla f_{{\cal N}_k}(x^k) $ and compute
$s^k$ satisfying (\ref{inexact1}).
\vskip 5 pt
\noindent
Step 3.  Find the smallest nonnegative integer $j$ such that
(\ref{lsnm}) holds for  $ t_k=2^{-j} $\\
\hspace*{33pt}     and set $ x^{k+1} = x^{k} + t_k s^k. $
}\label{GIN-R}

%
%
%
%

\begin{teo} \label{sr1}
Assume that A1 holds and  let $ \{x^k\} $ be generated by Algorithm GIN-R. If  $ \{\nu_k\} $ converges to zero  R-linearly then  $ \{x^k\} $ converges R -linearly  to $ x^*. $
\end{teo}




Next we show another important intermediate result.

\begin{lem} \label{etagama61}
Suppose that the assumptions of Theorem \ref{sr1} are satisfied and let $ \{x^k \} $ be a sequence generated by Algorithm GIN-R.  If $\eta_k$ is defined by (\ref{eta2}) then there exist positive constants $B_{1}$ and $ B_{2}$ and $\tau \in (0,1)$ such that $$\eta_{k} \leq B_{1} \tau^k \quad \mbox{and} \quad \gamma_{k} \leq B_{2} \tau^k$$ for all $k$ large enough. If $\eta_k=\bar{\eta}$ then $\gamma_{k}=C \bar{\eta}$ for all $k$ large enough.
\end{lem}

{\em Proof.} First, Theorem \ref{sr1} implies that the sequence of iterates $x^{k}$ converges to the unique solutions R-linearly -   for all $k$ large enough we have
\begin{equation} \label{err}
\|x^{k}-x^{*}\| \leq B \tau^k
\end{equation}
where   $B>0$ and $\tau \in (0,1)$.   Now, using the Taylor's expansion, Lemma \ref{Lemma3a} and (\ref{convexity}) we obtain for some $\theta_{k} \in  [x^{k-1}, x^k]$
\begin{eqnarray}
&  & |f_{\cal N}(x^{k})-m_{k-1} (s^{k-1})| \nonumber \\
&= & |f_{\cal N}(x^{k-1})+t_{k}\nabla f_{\cal N}(x^{k-1})^T s^{k-1}+\frac{1}{2} t_{k}^2(s^{k-1})^T \nabla^2 f_{\cal N}(\theta_k) s^{k-1} \nonumber \\
&-& f_{\cal N}(x^{k-1}) - \nabla f_{\cal N}(x^{k-1})^T s^{k-1} - \frac{1}{2} (s^{k-1})^T \nabla^2 f_{{\cal D}_{k-1}}(x^{k-1}) s^{k-1}| \nonumber \\
&\leq &  \frac{1-t_k}{\lambda_{1}}\|\nabla f_{\cal N}(x^{k-1})\|^2+\lambda_n \|s^{k-1}\|^2\nonumber \\
&\leq &  \frac{1}{\lambda_{1}}\|\nabla f_{\cal N}(x^{k-1})\|^2+\frac{\lambda_{n}}{\lambda^2_{1}} \|\nabla f_{\cal N}(x^{k-1})\|^2 \nonumber
\end{eqnarray}
Therefore, using (\ref{Lemma2}), for $\eta_{k}$ given by (\ref{eta2}) we obtain
\be \label{65} \eta_{k} \leq \frac{|f_{N}(x^{k})-m_{k-1} (s^{k-1})|}{\|\nabla f_{N}(x^{k-1})\|}
\leq (\frac{1}{\lambda_{1}}+\frac{\lambda_{n}}{\lambda^2_{1}})\lambda_{n} \|x^{k-1}-x^{*}\| \leq B_{1} \tau^{k} \ee
where $B_1= (\frac{1}{\lambda_{1}}+\frac{\lambda_{n}}{\lambda^2_{1}})\lambda_{n}B /\tau$.
Moreover, for $\gamma_k$ given by (\ref{62}) we get
$$\gamma_{k} \leq C \max \{B_{1} \tau^{k},\lambda_{n} B  \tau^{k}\}=C \max \{B_{1} ,\lambda_{n} B \}\tau^{k} :=B_2 \tau^{k}.$$
Considering the case with $\eta_{k}=\bar{\eta}$,  since the gradient converges to zero, for  all $k$ large enough we have $\|\nabla f_{\cal N}(x^{k})\| < \bar{\eta}$ and therefore
$\gamma_{k}=C \bar{\eta}$.
$\square$

Now, we are ready to show the main result of this subsection.



\begin{teo} \label{mss-super}
Suppose that the assumptions of Theorem \ref{sr1} and Assumption A2 are satisfied.
Moreover, assume $c\in (0,1/4)$ in (\ref{lsnm}) and $0 < C <(1/{\bar{\eta}} - 1)\lambda_1/2$ in (\ref{Dk_choice}) .
 Let $ \{x^k\} $ be a sequence generated by Algorithm GIN-R.
Then, there are positive constants $V_{1}, V_{2}, C_1, C_2$ and $\tau \in (0,1)$ such that for all $k$ sufficiently large
\begin{itemize}
\item[a)] if  $\eta_k$ is defined by (\ref{eta2}) then  $$E(\|x^{k+1}-x^{*}\|^2) \leq \left(V_{1} \tau^{2k}+V_{2}  \alpha_{k}\right) E(\|x^{k}-x^{*}\|^2);$$
\item[b)] if  $\eta_k=\bar{\eta}$ is sufficiently small then $$E(\|x^{k+1}-x^{*}\|^2) \leq \left((C_{1} \tau^{k}+C_{2} \bar{\eta})^{2}+V_{2} \alpha_{k}\right) E(\|x^{k}-x^{*}\|^2).$$
\end{itemize}
\end{teo}

{\em Proof.}
 Since the assumptions of Theorem \ref{sr1} are satisfied, the sequence $\{x^k\}$ converges to the solution R-linearly, independently of $D_{k}$ and for $\eta_k \in (0,1)$. Thus, there exist constants $B>0$ and $\tau \in (0,1)$ such that $\|x^{k}-x^{*}\| \leq B \tau^k$. Moreover, $N_{k}=N$ for all $k$ sufficiently large. So, without loss of generality, we assume that the full sample is used for the gradient and the function and $\|\nabla f_{\cal N}(x^k)\|<\bar \eta$.

Employing (\ref{conv1}) and Lemma \ref{Lemma3a}, we have the following estimate
\begin{eqnarray} \label{642}
 \|x^{k}+t_k s^{k}-x^{*}\|^2
&\leq  & \frac{2}{\lambda_{1}} (f_{\cal N}(x^{k}+t_k s^{k})-f_{\cal N}(x^*)) \nonumber \\
& \leq & \frac{2}{\lambda_{1}} (f_{\cal N}(x^k)-f_{\cal N}(x^*)+t_k (\nabla f_{\cal N}(x^k))^T s^{k}+\frac{\lambda_{n}}{2} \|t_k s^k\|^2) \nonumber \\
&\leq  & \frac{2}{\lambda_{1}} (\frac{1}{\lambda_{1}}\|\nabla f_{\cal N}(x^k)\|^2+\frac{\lambda_{n}}{2} \frac{1}{\lambda^2_{1}} \|\nabla f_{{\cal N}_k} (x^k)\|^{2}) \nonumber \\
&\leq  & \frac{2}{\lambda_{1}} (\frac{1}{\lambda_{1}}+ \frac{\lambda_{n}}{2 \lambda^2 _{1}})\lambda_{n}^{2}\|x^{k}-x^{*}\|^{2} :=
 V_{2}  \|x^{k}-x^{*}\|^{2}.
\end{eqnarray}

Let us denote by $A_{k}$ the event $\|\nabla^2 f_{{\cal D}_{k}}(x^k)-\nabla^2 f_{\cal N}(x^k)\|\leq \gamma_{k} .$
Due to  Step 1.c of algorithm GIN-R
it follows  that $P(A_{k}) \geq 1-\alpha_{k}$, i.e., $P(\bar{A}_{k}) \leq \alpha_{k}$. Notice that (\ref{642}) holds in both cases but in the case of $ A_k $ we can derive better estimate.

Assume that $A_{k}$ happens. Then
\begin{eqnarray} \label{641}
&  &  \|\nabla^2 f_{\cal N}(x^k)s^k+\nabla f_{\cal N}(x^k)\| \nonumber \\
& \leq & \|(\nabla^2 f_{\cal N}(x^k)-\nabla^2 f_{{\cal D}_{k}}(x^k))s^k\|+\|\nabla^2 f_{{\cal D}_{k}}(x^k)s^k+\nabla f_{\cal N}(x^k)\|\nonumber \\
&\leq  &  \gamma_{k} \|s^{k}\|+ \eta_{k}\| \nabla f_{\cal N}(x^k)\| \leq (\gamma_{k} /\lambda_{1}+ \eta_{k})\| \nabla f_{\cal N}(x^k)\|\nonumber.
\end{eqnarray}
Moreover, as $\gamma_k< C\bar \eta$ we can repeat the reasoning used in the proof of Theorem \ref{thRrate} to conclude that the full step is accepted for $k$ sufficiently large. Also, notice that the standard assumptions for (Inexact) Newton method are satisfied, so following the standard steps for IN analysis one can prove
\be \label{novo1} \|x^{k}+s^{k}-x^{*}\| \leq c_{1} (\|x^{k}-x^{*}\|+\gamma_{k} /\lambda_{1}+ \eta_{k})\|x^{k}-x^{*}\|\ee
for some positive constant $c_{1}$ and for $k$ sufficiently large. 
Indeed, denoting $\tilde{r}_k=\nabla^2 f_{\cal N}(x^k)s^k+\nabla f_{\cal N}(x^k)$, it holds  $\|\tilde{r}_k\| \leq (\gamma_{k} /\lambda_{1}+ \eta_{k})\| \nabla f_{\cal N}(x^k)\|$ and
 $$\|x^{k}+s^{k}-x^{*}\|=\|x^{k}+(\nabla^2 f_{\cal N}(x^k))^{-1} \tilde{r}_k-(\nabla^2 f_{\cal N}(x^k))^{-1}\nabla f_{\cal N}(x^k) -x^{*}\|.$$ Since the Newton's method converges quadratically there exists $\kappa>0$ such that
 $$\|x^{k}-(\nabla^2 f_{\cal N}(x^k))^{-1}\nabla f_{\cal N}(x^k) -x^{*}\|\leq \kappa \|x^k-x^*\|^2.$$
 Now, using $\|(\nabla^2 f_{\cal N}(x^k))^{-1}\| \leq 1/\lambda_1$,  $\|\nabla f_{\cal N}(x^k)\| \leq  \lambda_n \|x^{k}-x^{*}\| $ and defining $c_1=\max \{\kappa, \lambda_n/\lambda_1\}$ we obtain (\ref{novo1}). 
Using   inequality  (\ref{err}) and squaring inequality (\ref{novo1}) we obtain
\be \label{6ak}
\|x^{k}+s^{k}-x^{*}\|^2 \leq 2 c_{1}^{2}  (B^2 \tau^{2 k}+  (\gamma_{k} /\lambda_{1}+ \eta_{k})^2) \|x^{k}-x^{*}\|^2.
\ee

Now, we distinguish two cases depending on $ \eta_k. $ Using the result of Lemma \ref{etagama61} and assuming that $k$ is sufficiently large we obtain the following.
\begin{itemize}
\item[a)] If $\eta_k$ is defined by (\ref{eta2}), for   $V_1=2 c_{1}^2 (B^2+2 B_1^2+2 (B_{2}/\lambda_1)^2)$ we get
\begin{eqnarray} \label{643}
&  &  E(\|x^{k+1}-x^{*}\|^2) \nonumber \\
&= & P(A_{k}) E(\|x^{k+1}-x^{*}\|^2| A_{k} )+P(\bar{A}_{k}) E(\|x^{k+1}-x^{*}\|^2| \bar{A}_{k} )\nonumber \\
&\leq  & (V_{1} \tau^{2k} +\alpha_{k}V_{2}) E(  \|x^{k}-x^{*}\|^{2}).
\end{eqnarray}
\item[b)] Considering $\eta_k=\bar{\eta}$, for  $C_{1}=2 c_{1}^{2}  B^2$ and $ C_2=2 c_{1}^{2} (C  /\lambda_{1}+ 1)^2$ we get
\be \label{6243} E(\|x^{k+1}-x^{*}\|^2) \leq \left(C_{1} \tau^{2k}+C_{2}  \bar{\eta}^{2}+V_{2} \alpha_{k}\right) E(\|x^{k}-x^{*}\|^{2}).\ee
\end{itemize}
$\square$

We conclude the analysis with the following corollary.
\begin{cor} \label{cor-mss-super}
Assume that the conditions of Theorem \ref{mss-super} hold and let $ \{x^k\} $ be a sequence generated with Algorithm GIN-R. Then the sequence $\{x^k\} $ converges to $ x^* $ in the mean square sense:
\begin{itemize}
\item[a)] linearly  if $\eta_k=\bar{\eta}$ is sufficiently small and $\alpha_{k}< \frac{1-C_{2}^2 \bar{\eta}^2}{V_{2}};$
\item[b)] superlinearly  if  $\eta_k$  is defined by (\ref{eta2}) and
$\lim_{k \rightarrow \infty} \alpha_{k}=0.$
\end{itemize}
\end{cor}

Note that the requirement on the probability $\alpha_k$ in order to get the $q$-linear convergence   is not too demanding  in what concerns the sample size,
because it influences only the logarithmic factor   in \eqref{minbound}.

\subsection{Unbounded sample - GIN method}

In many applications, the number of training points is enlarged over time so the cardinality of the sample set  $N$ is actually unbounded. This motivated us to consider the following problem as well
\be \label{ini1}
\min_{x \in \mathbb{R}^n} f(x)=E(F(x,\xi)),
\ee
where $\xi$ is a random variable defined on a probability space $({\cal A}, {\cal F}, P)$ and $F$ is twice differentiable function with respect to $x$. Let us  denote
 $f_{i}(x):=F(x,\xi_i)$ where $\xi_i, \; i=1,2,...$ is an i.i.d. sequence of variables following the same distribution as $\xi$. For example, $\xi_i$ can represent the pair of input-output variables  in machine learning problems.
Then, we can use the same notation as in the previous sections to define SAA approximation of the objective function and its derivatives: $f_{{\cal N}}, \nabla f_{{\cal N}}, \nabla^2 f_{{\cal D}}$. We will prove that, under appropriate assumptions, GIN  converges in m.s. towards the solution of problem (\ref{ini1}).

Let us formalize the assumption about random variables $\xi_i$.

\noindent {\bf{Assumption B1}} $\xi_i, \; i=1,2,...$ is an i.i.d. sequence of variables.

Next we assume that the sequence of iterates $\{x^k\}$  belongs to a bounded set. This assumption might sound like a strong one, but it is needed for convergence theory and  the analysis of properties of $ F(x,\xi) $ that guaranty this property is beyond the scope of this paper. Notice that a weaker assumption is used in \cite{bbn} (the assumption of bounded moments of iterates). However, here we are using a linesearch procedure to compute the steplength $t_k$ and therefore we need to use the SAA approximation for the function. On the contrary,   \cite{bbn} emploies  a fixed steplength whose definition require the knowledge of the maximum eigenvalue $\lambda_n$, but allows to get rid of computing the objective function.

\noindent {\bf{Assumption B2}} There exists a compact set $\Omega$ such that $\{x^k\}_{k \in \mathbb{N}} \subseteq \Omega$.

\noindent {\bf{Assumption B3}} $F(\cdot, \xi) \in C^2(\mathbb{R}^n)$ for every $\xi$. $F$ and $\nabla F$ are dominated by an integrable functions $M_f(\xi)$ and $M_g(\xi)$, respectively, on an open set containing $\Omega$.

Assumption B1 implies that $E(f_{{\cal N}}(x))=f(x)$. Moreover,  B1 and B3 imply that $\nabla f(x)=E(\nabla F(x, \xi))$ and therefore $E(\nabla f_{{\cal N}}(x))=\nabla f(x)$. Furthermore, the Uniform law of large numbers implies that $ f_{{\cal N}}$ and $\nabla f_{{\cal N}}$ almost surely (a.s.) converge to $f(x)$ and $\nabla f(x)$, respectively, uniformly on $\Omega$ when $N$ tends to infinity. Let us write this down more formally. Denote
\be \label{inie} \ef = \max_{x \in \Omega} |f_{{\cal N}}(x)-f(x)|, \quad \eg = \max_{x \in \Omega} \|\nabla f_{{\cal N}}(x)-\nabla f(x)\|. \ee
 Then, $\lim_{N \rightarrow \infty} \ef=0$ and $\lim_{N \rightarrow \infty} \eg =0$, a.s. and using the Lebesgue Dominated Convergence Theorem  (see Theorem 7.31 of \cite{Sapiro}) we obtain
\be \label{ini2}
\lim_{N \rightarrow \infty } E(\ef)=0, \quad \lim_{N \rightarrow \infty } E(\eg)=0.
\ee

Assuming the strong convexity of $f_i$ as in assumption A1, it is easy to show that $f_{{\cal N}}$ is also strongly convex with the same constants $\lambda_1$ and $\lambda_n$ for any ${\cal N}$. Moreover, assuming B3, $f$ also remains strongly convex with the constant $\lambda_1$. Indeed, for an arbitrary $i$ and $x, y$ there holds
$$f_i (y)\geq f_i (x)+ \nabla^T f_i (x) (y-x)+\frac{\lambda_1}{2} \|x-y\|^2.$$
Taking the expectation and using that $E(\nabla f_i (x))=\nabla f(x)$ we obtain the strong convexity of $f$. Therefore, problem (\ref{ini1}) has an unique solution $x^*$.

\begin{teo} \label{ini3}
Suppose that the assumptions A1, B1-B3 hold and that  $N_{k} \to \infty$.  Then any  sequence $\{x^k\}_{k \in \mathbb{N}}$ generated by GIN converges towards the solution of the problem (\ref{ini1}) in the mean square  sense.
\end{teo}

{\em Proof.}
First, notice that assumptions B1-B3 imply that $|f_{{\cal N}}(x^k)-f(x^k)| \leq \ef$ for every $k$, so following the reasoning as in the proof of Theorem \ref{thglobal} we obtain
\be \label{ini6} f(x^{k+1}) - f(x^*) \leq \fnk(x^k)-c \bar{t} q \|\nabla \fnk (x^k)\|^2+\efk- f(x^*)+\nu_k,\ee
where $\bar{t}=(1-c)\lambda_1/\lambda_n $ and $q= (\lambda_1 (1-\bar{\eta})^2)/(\lambda^2_{n})$.
Now,
\begin{eqnarray}
\|\nabla \fnk (x^k)\|^2 & = & \|\nabla \fnk (x^k)-\nabla f(x^k) + \nabla f(x^k)\|^2 \nonumber \\
&=& \| \nabla f(x^k)\|^2+2 (\nabla f(x^k))^T (\nabla \fnk (x^k)-\nabla f(x^k)) \nonumber \\
& + & \|\nabla \fnk (x^k)-\nabla f(x^k)\|^2\nonumber\\
&\geq & \| \nabla f(x^k)\|^2-2 \|\nabla f(x^k)\| \| \nabla \fnk (x^k)-\nabla f(x^k)\| \nonumber\\
&\geq & \| \nabla f(x^k)\|^2-2 M_g \egk,
\end{eqnarray}
where the last inequality follows from (\ref{inie}), continuity of $\nabla f$,  Assumption B2 and $M_g=\max_{x \in \Omega} \|\nabla f(x)\|$.
On the other hand, strong convexity of $f$ implies that
$-\| \nabla f(x^k)\|^2\leq - \lambda_1 (f(x^k)-f(x^*))$. Putting all together  into (\ref{ini6}) we obtain
$$f(x^{k+1}) - f(x^*) \leq (f(x^{k}) - f(x^*)) (1-\omega) + 2\efk+2 c \bar{t} q M_g \egk+\nu_k,$$
where $\omega = c \bar{t} q \lambda_1 \in (0,1)$.  Applying expectation we get
$$E(f(x^{k+1}) - f(x^*) )\leq E(f(x^{k}) - f(x^*)) (1-\omega) + a_k,$$
where $a_k= 2E(\efk)+2 c \bar{t} q M_g E(\egk)+\nu_k$. Now, (\ref{ini2}), (\ref{errorsequence}) and the assumption that $N_{k} \to \infty$  together imply that $\lim_{k \rightarrow \infty} a_k =0$. Therefore, it follows (see \cite{dj}) that
$$\lim_{k \rightarrow \infty} E(f(x^{k}) - f(x^*))=0. $$
Finally, strong convexity implies $\|x^k-x^*\|^2 \leq (f(x^{k}) - f(x^*)) 2/\lambda_1$ thus
$$\lim_{k \rightarrow \infty} E(\|x^k-x^*\|^2)=0 $$
$\square$

\subsection{Relaxing the strong convexity - method GIN-M}

In this part of the paper, we consider a relaxation of the strong convexity assumption (\ref{convexity}) by letting $\nabla^2 f_i(x)$ be only positive semidefinite while the final objective function remains strongly convex. Similar assumptions are stated in  \cite{mit1}. 
Notice that Theorem \ref{ini3} does not impose any assumption on the size of the Hessian subsample. But relaxation of the strong convexity imposes condition on the Hessian sample size  needed to obtain a  positive definite Hessian in probability. We  use the bound provided in \cite{mit1} but  the approach presented here differs in several ways from the one in \cite{mit1}.  First of all, we  prove mean square convergence to the solution of (\ref{ini1}) under appropriate conditions given in the sequel, while in \cite{mit1} convergence with some (high) probability for finite sum problems like (\ref{problem1}) is considered. Second, we continue to use  CG as the inner solver - although modified in this case to cope with possibly singular matrix.

\noindent {\bf{Assumption C1}}  The functions $f$ and $ F(\cdot, \xi)$ are twice continuously differentiable and there exist $0<\lambda_1<1$ and $\lambda_n>0$ such that for every $x, \xi$
$$0  \preceq \nabla_{x}^2 F(x, \xi) \preceq \lambda_n I \quad \mbox{and} \quad \lambda_1 I \preceq \nabla^2 f(x) \preceq \lambda_n I.$$

This assumption ensures that the unique solution of the original problem still exists. Moreover, we assume that the Hessian approximations are unbiased.

\noindent {\bf{Assumption C2}}  For every $x \in \mathbb{R}^n$ and every $i\in \mathbb{N}$ there holds $E(\nabla^2 f_i(x))=\nabla^2 f(x). $

This assumption allows us to use the Matrix Chernoff result (see \cite{tropp12}) and to obtain the bound presented in Lemma 1 of \cite{mit1}. Although we observe unbounded sample, the same result holds. More precisely, under the assumptions C1 and C2 we obtain that, given $ \mu \in (0,1-\lambda_1)$ the following holds

$$P\left(\lambda_{\min}(\nabla^2 f_{{\cal D}}(x))\geq \mu \right)\geq 1-\alpha$$
if
\be \label{ini12} D \geq \frac{2 \lambda_n (1-\lambda_1)^2  \ln(n/\alpha)}{\mu^2 \lambda_1}:=\bar{D}(\alpha).\ee
Given that the subsampled Hessian $ \nabla^2 f_{{\cal D}_k} $ might be singular, for the computation of the step $s^k$ at Step 2 of Algorithm GIN-R,  we proceed as follows.
If at iteration $j$ of CG it happens that  $(s^k_j)^T \nabla^2 f_{{\cal D}_{k}}(x^k) s^k_j=0$ we stop CG and set $s^k=s^k_{j-1}$. Notice that $ s^k $ is still a descent direction as $
(s^k_{j-1})^T\nabla f=- (s^k_{j-1})^T \nabla^2 f_{{\cal D}_{k}}(x^k) s^k_{j-1}<0$.
The algorithm we use is again GIN-R  where  the sample size is selected such that subsampled Hessian is positive definite with probability $1-\alpha. $ The modified CG, as  explained above, is used. We refer to the obtained procedure as Algorithm GINR-M and sketch its iteration $k$ in Algorithm \ref{GIN-M}.
{\bf
 \algo{GINR-M}{$k$-th iteration of Method GINR-M}{
\noindent
Given  $ x^k \in {\mathbb R}^n, \;  c\in (0,1), \bar{\eta} \in (0,1), \; C>0, \; \{\nu_k\},\; \alpha\in (0,1)\; \mu\in (0,1-\lambda_1).  $
\begin{description}
\item{\bf Step 1.a} Choose $ {\cal N}_k, \eta_k\in (0,\bar \eta). $
\item{\bf Step 1.b}  If $N_k\leq \bar{D}(\alpha)$ given in \eqref{ini12},
set ${\cal D}_k = {\cal N}_k$. Else, choose  {\bf  $D_k\ge  \bar{D}(\alpha)$} and the  subsample $ {\cal D}_k $ randomly and uniformly from  $ {\cal N}_k  $.
\item{\bf Step 2. } Determine $ s^k $ with modified CG: if $(s^k_j)^T \nabla^2 f_{{\cal D}_{k}}(x^k) s^k_j=0$  for some inner iteration $j$, set  $s^k=s^k_{j-1}$. Otherwise, find the step $ s^k $ such that (\ref{inexact1}) holds.
\item{\bf Step 3. } Find the smallest nonnegative integer $j$ such that (\ref{lsnm}) holds for  $ t_k=2^{-j} $  and set $ x^{k+1} = x^{k} + t_k s^k. $
\end{description}
}\label{GIN-M}
}

Relaxing the strong convexity results in losing the usual relation between the step and the gradient stated in Lemma \ref{Lemma3a}.
Therefore, we need the following assumption.

{\bf{Assumption C3}} There exists a constant $M_s>0$ such that the step generated by GINR-M satisfies  $\|s^k\|\leq M_s$ for every $k$.

 A comment is due with respect to the above Assumption.  Assume that  the step $s^k$ computed at Step 4 of GINR-M has been generated at iteration $j$ of CG. Then it belongs to the Krylov subspace
$${\cal K}_j=span\{\nabla f_{{\cal N}_k}(x^k), (\nabla^2 f_{{\cal D}_{k}}(x^k))\nabla f_{{\cal N}_k}(x^k),\ldots,
(\nabla^2 f_{{\cal D}_{k}}(x^k))^{j-1}\nabla f_{{\cal N}_k}(x^k) \}.$$ Let $V$ be an orthonormal basis of ${\cal K}_j$, then
$s^k=Vy$,  where $y\in \mathbb{R}^j$. Therefore
$$
(s^k)^T \nabla^2 f_{{\cal D}_{k}}(x^k) s^k=y^T V^T\nabla^2 f_{{\cal D}_{k}}(x^k) V y \ge  \lambda_{\min} (V^T \nabla^2 f_{{\cal D}_{k}}(x^k) V)  \|y\|^2.
$$
Then, noting that $\|s\|=\|y\|$ as $V$ is orthornormal,
from
$$
(s^k)^T \nabla^2 f_{{\cal D}_{k}}(x^k) s^k=-(s^k)^T\nabla f_{{\cal N}_k}(x^k)
$$
and the Cauchy-Schwartz inequality, there follows
$$
\|s^k\|\le \frac{1}{ \lambda_{\min} (V^T \nabla^2 f_{{\cal D}_{k}}(x^k) V) }  \|\nabla f_{{\cal N}_k}(x^k)\|.
$$
So, Assumption C3 is satisfied whenever the minimum eigenvalue of the projected subsampled Hessian   $V^T \nabla^2 f_{{\cal D
}_{k}}(x^k) V$ is bounded away from zero and Assumption B2 holds.

\begin{teo} \label{ini13}
Suppose that the assumptions C1-C3, B1-B3 hold and that  $N_{k}$ tends to infinity. Then   there exist $\alpha$ small enough
 such that any sequence $\{x^k\}_{k \in \mathbb{N}}$ generated by GINR-M converges towards the solution of  (\ref{ini1}) in m.s.
\end{teo}

{\em Proof.}
Since $N_{k}$ tends to infinity, $N_k\geq \bar{D}(\alpha)$ will be satisfied for all $k$ large enough ($k \geq k(\alpha)$) and thus in Step 1.b $D_k$ is chosen such that (\ref{ini12}) holds. Since we are interested in asymptotic result, without loss of generality we assume that $k \geq k(\alpha)$.
Denote by $A_k$ an event $\lambda_{\min}(\nabla^2 f_{{\cal D}_{k}}(x^k))\geq \mu$. Since  $D_k$ satisfies (\ref{ini12}) it follows that  $P(\bar{A}_k)\leq \alpha$.


Assume that $A_k$ happens. Then, we can  proceed as in the proof of Theorem \ref{ini3} as
$f$ is strongly convex and   $\lambda_{\min}(\nabla^2 f_{{\cal D}}(x^k))\geq \mu$ yields  $(\nabla\fnk (x^k))^T s^k \leq -\mu \|s^k\|^2$. We obtain
$$f(x^{k+1}) - f(x^*) \leq (f(x^{k}) - f(x^*)) (1-\omega) + 2\efk+\theta \egk+\nu_k,$$
where  $\omega = c \bar{t} q \lambda_1 \in (0,1)$, $\bar{t}=(1-c)\mu /\lambda_n$, $q=\mu (\frac{1-\bar{\eta}}{\lambda_n})^2$,  $\theta=2 c \bar{t} q M_g$ and $M_g=\max_{x \in \Omega} \|\nabla f(x)\|$.


On the other hand, assume that $\bar{A}_k$ happens. Then by the Taylor expansion and assumption C1 we obtain
\be \label{cbar} f(x^{k+1})-f(x^*) \leq f(x^{k})-f(x^*)+t_k (\nabla f (x^k))^T  s^k+\frac{1}{2} \lambda_{n} \|t_k s^k\|^2.\ee
Again, using the strong convexity of $f$  we get
\begin{eqnarray}  \label{ini14} \|t_k s^k\|^2 & = & \|x^{k+1}-x^k\|^2\leq 2 (\|x^{k+1}-x^*\|^2+\|x^{k}-x^*\|^2) \nonumber \\
& \leq &  \frac{4}{\lambda_1}(f(x^{k+1})-f(x^*)+f(x^{k})-f(x^*)).
\end{eqnarray}
Moreover,
\begin{eqnarray}
(\nabla f (x^k))^T  s^k & = &   (s^k)^T ( \nabla  f (x^k) \pm \nabla \fnk (x^k)) \nonumber \\
&\leq & \| s^k\| \|\nabla  f (x^k) - \nabla \fnk (x^k)\|+(\nabla \fnk (x^k))^T s^k \nonumber\\
&= & \| s^k\| \|\nabla  f (x^k) - \nabla \fnk (x^k)\|- (s^k)^T \nabla^2 f_{{\cal D}_{k}}(x^k)s^k \nonumber\\
&\leq & M_s \egk, \label{Ms}
\end{eqnarray}
where the last inequality follows from the fact that $\nabla^2 f_{{\cal D}_{k}}(x^k)$ is positive semidefinite.
Putting (\ref{Ms})  into (\ref{cbar}),  together with (\ref{ini14}) we obtain
$$f(x^{k+1})-f(x^*) \leq (f(x^{k})-f(x^*))(1+2 \lambda_n / \lambda_1)+2 \lambda_n / \lambda_1 (f(x^{k+1})-f(x^*))+M_s \egk.$$

Combining all together we get
\begin{eqnarray*}
&  & E(f(x^{k+1}) - f(x^*))   \\
& =  &   P(A_{k}) E(f(x^{k+1}) - f(x^*)| A_{k} )+P(\bar{A}_{k}) E(f(x^{k+1}) - f(x^*)| \bar{A}_{k} )  \\
&\leq & E(f(x^{k}) - f(x^*)) (1-\omega) + 2E(\efk)+\theta E(\egk)+\nu_k \\
&+ &\alpha \left(E(f(x^{k})-f(x^*))(1+2 \lambda_n / \lambda_1)+2 \lambda_n / \lambda_1 E(f(x^{k+1})-f(x^*)) \right)  \\
& + &  \alpha M_s E(\egk).
\end{eqnarray*}
Rearranging the previous inequality and assuming $\alpha < \frac{\lambda_1} {2 \lambda_n}$ we obtain
$$E(f(x^{k+1}) - f(x^*)) \leq \tau E(f(x^{k}) - f(x^*))+a_k,$$
where $$\tau=\frac{1-\omega +\alpha (1+2 \lambda_n/\lambda_1)}{u}, \quad u=1-\alpha 2 \lambda_n / \lambda_1$$
and
$$a_k=\frac{1}{u} (\nu_k+2E(\efk)+(\theta+ \alpha M_s)  E(\egk)). $$
Notice that $\tau \in (0,1)$ provided that  $\alpha$ is small enough. Moreover, as discussed in the previous proof, $a_k$ tends to zero and the result follows.
$\square$

\section{Numerical results}

In this section we report on our numerical experience with subsampled inexact Newton approaches. Our experiments were performed in Matlab R2017a, on a
Intel Core i5-6600K CPU 3.50 GHz x 4 16GB RAM.
For the approximate solution of the arising linear systems we used the conjugate gradient method (CG) implemented in  the Matlab function {\tt pcg}. We did not employ a preconditioner and we used the conjugate gradient method in a matrix-free manner. Then, only the product of $\nabla f^2_{  {   \cal D} _k}$ with vectors  is needed.
The aim of this section it to provide numerical evidence of the benefits deriving by the employment of adaptive rules, streaming out from our theory, for  choosing forcing terms and Hessian sample size.Then,
in our numerical experimentation we use full gradients and functions, i.e. $N_k=N$ for $k>0$ and we  compare the full Inexact Newton ({\tt FIN}) method with $\eta=10^{-4}$, the subsampled Hessian ({\tt SIN}) method with $\eta_k=10^{-4}$ and $D_k=0.3 N$  for all $k$'s, the subsampled inexact method with   adaptive choices  of $\eta_k$'s
and $D_k=0.3 N$  for all $k$'s ({\tt SINA\_FT}) and   the subsampled inexact method with  adaptive choices  of $\eta_k$'s  and $D_k$'s ({\tt SINA\_FT\_Dk}).
We also consider a subsampled method with constant $D_k=0.3N$ and a maximum number  of five iterations allowed to {\tt pcg} ({\tt SIN\_cg5}).

In  {\tt SINA\_FT} and {\tt SINA\_FT\_Dk},   $\eta_k$ is  chosen as follows:
$$
\left \{
\begin{array}{l}
 \eta_k=\min\{0.1,\max\{|f(x_k)-m_{k-1}(s^k)|/\|\nabla f(x_{ k-1})\|\},10^{-3})\}\\
\eta_0=0.1.
\end{array}
\right .
$$
with $ m_{k-1}(s^k) $ given in \eqref{genmod}
This choice  is made   according to  (\ref{eta2}).
Finally, we choose  the sample size $D_k$ in  {\tt SINA\_FT\_Dk} as:
$$
D_k=\lceil{\max\{c_0 D_0, \min\{c_1  \min\{\frac{1}{\eta_k^2},\frac{1}{\|\nabla f(x_k)\|^2}\},N\}}\rceil \;\;\;c_0,c_1>0
$$
The above rule is based on inequalities   (\ref{minbound})  and (\ref{62}). The idea is to choose $D_k$ inversely proportional to $\gamma_k^2$ given in (\ref{62}) as suggested by the bounds  in (\ref{minbound}). The choice of costants $c_0$ and $c_1$ depends on the convergence behaviour of CG method. In fact, if CG converges fast large values of $D_k$ should be used, as the lost in the convergence rate due to less accurate second order information is not compensated by the reduced cost of Hessian-vector products. On the other hand, when CG is slower, smaller values of $D_k$ must be used as a large number of matrix-vector products are needed.
Then, $c_0$ and $c_1$ are chosen according to the following strategy.  We set $c_0=1$, $c_1=0.05$ in case at the previous iteration CG needed more than 20 iterations. Otherwise we set $c_0=2$ and $c_1=1$.
Moreover,  $D_0$ is set to $0.1 N$. This choice is motivated by the fact that we allow $D_k$ to change and increase, then we can start with a small $D_k$ leaving the method free to adaptively modify it.

In all the subsampled  methods  the set ${\cal D}_k$ is chosen randomly using the Matlab function   {\tt randperm}.
 All the methods under  comparison are in the framework of Algorithm GIN-R, i.e. the nonmonotone linesearch   (\ref{lsnm}), with $c=10^{-4}$ and $\nu_k=\max(1,f(x^0)/k^{1.1}$  is applied.

 The problem we consider is the binary classification problem. We suppose to have at disposal a training set composed of pairs $\{(a_i, b_i)\}$ with $a_i \in \IR^n$,
$b_i \in  \{-1, +1\}$ and $i = 1,\ldots,N$, where $b_i$ denotes the correct sample classification. We
perform a logistic regression, then we consider as a training objective function the logistic loss
with $\ell_2$ regularization, \cite{bbn}, i.e. in problem (\ref{problem1}) we have
\be \label{c}
f_i(x) = \log c(x,\xi_i)+\lambda \|x\|^2, \; c(x,\xi_i) =1+ e^{-b_i a_i^T x}
\ee
where
$ \xi_i = (a_i,b_i)$.   Furthermore, the gradients and the Hessians have special forms,
\begin{equation}\label{gradc}
\nabla f_i(x) = \frac{(1-c(x,\xi_i) ) }{c(x,\xi_i) } b_i a_i +2 \lambda  x
\end{equation}
\begin{equation}\label{Hc}
 \nabla^2 f_i(x) = -\frac{1-c(x,\xi_i)}{c^2(x,\xi_i)} a_i a_i^T+2\lambda I.
\end{equation}

Note that the evaluation of the full function $f$ requires the evaluation of the quantities
$(c(x,\xi_i) -1) b_i a_i$, for $i=1,\ldots,n$ and once these quantities have been computed they can be used for evaluating
$ \nabla f_i(x)$ and  $\nabla^2 f_i(x)$, $i=1,\ldots,n$. Then, due to the form of gradient and Hessian of each $f_i$ given in (\ref{gradc}) and (\ref{Hc}),  the evaluation of $\nabla f_i(x)$ comes from free and the evaluation of $\nabla^2 f_i$ times a vector is as expensive as evaluating $f_i(x)$. Then, we evaluate the performance of the methods under comparison in terms of     full function evaluations (FEV).
We underline that one  {\tt pcg} iteration costs as $\frac{D_k}{N}$ FEV.
Finally,  in (\ref{c}) we set $\lambda=1/N$.

We used the following three datasets 
\begin{itemize} 
\item {\tt CINA0}  \cite{CINA}, $N= 16033$, $n= 132$;
\item  {\tt Mushrooms} \cite{UCI}, $N=5000$, $n=112$;
\item   {\tt Gisette}  \cite{UCI},  $N=6000$,  $n=5000$.
\end{itemize}

We stopped the methods under comparison when
$$
||\nabla f(x_k)||\le 10^{-4}
$$
or when a maximum number of 50 nonlinear iterations is reached. In order to compute the training error we compute the minimizer $x^*$ with the full Newton method to a tight accuracy, i.e we run it until  $||\nabla f(x_k)||$ is less than $10^{-8}$.

\begin{figure}
\hspace*{-28pt}\includegraphics[height=0.25\textheight, width=0.6\textwidth]{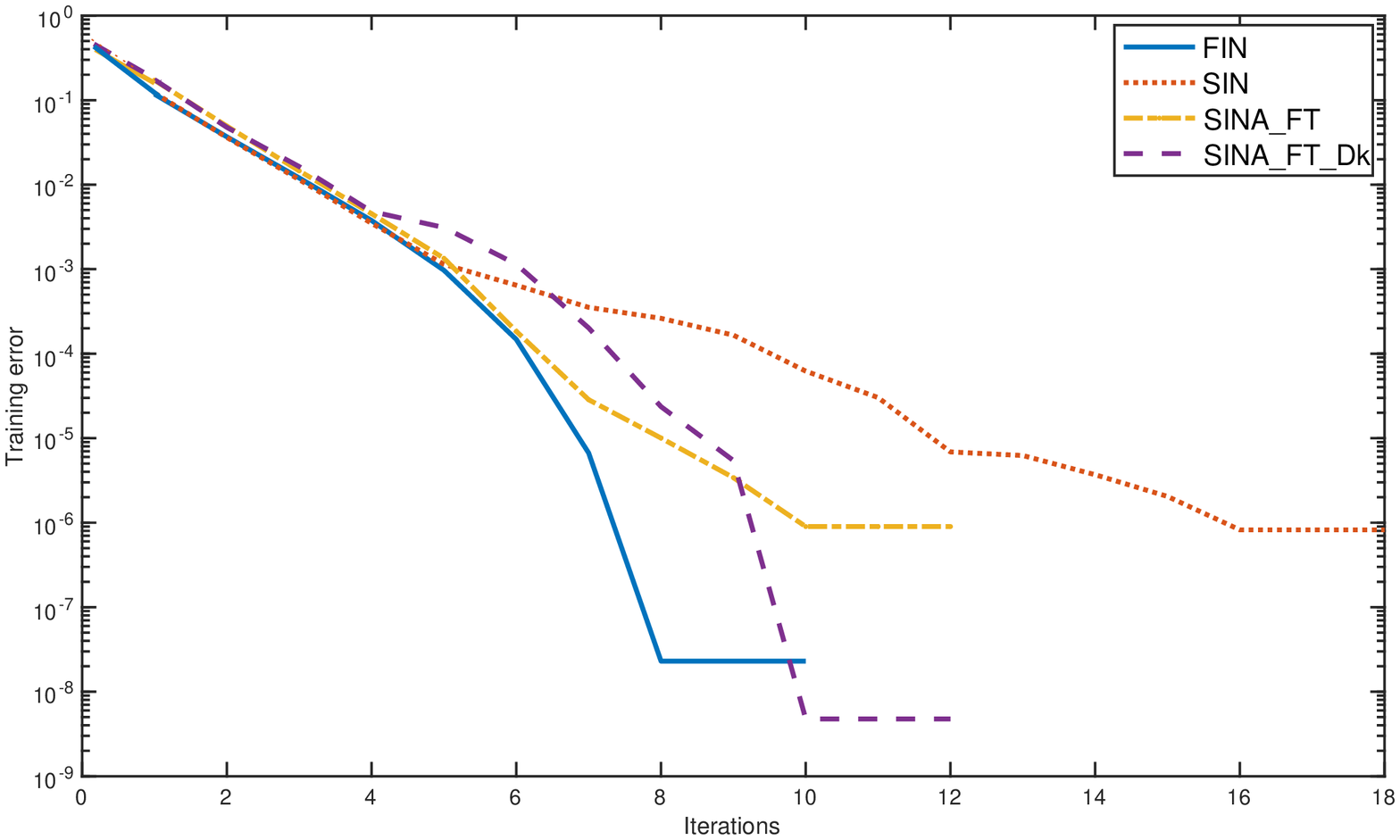}
\hspace*{-15pt}\includegraphics[height=5.2cm, width=0.6\textwidth]{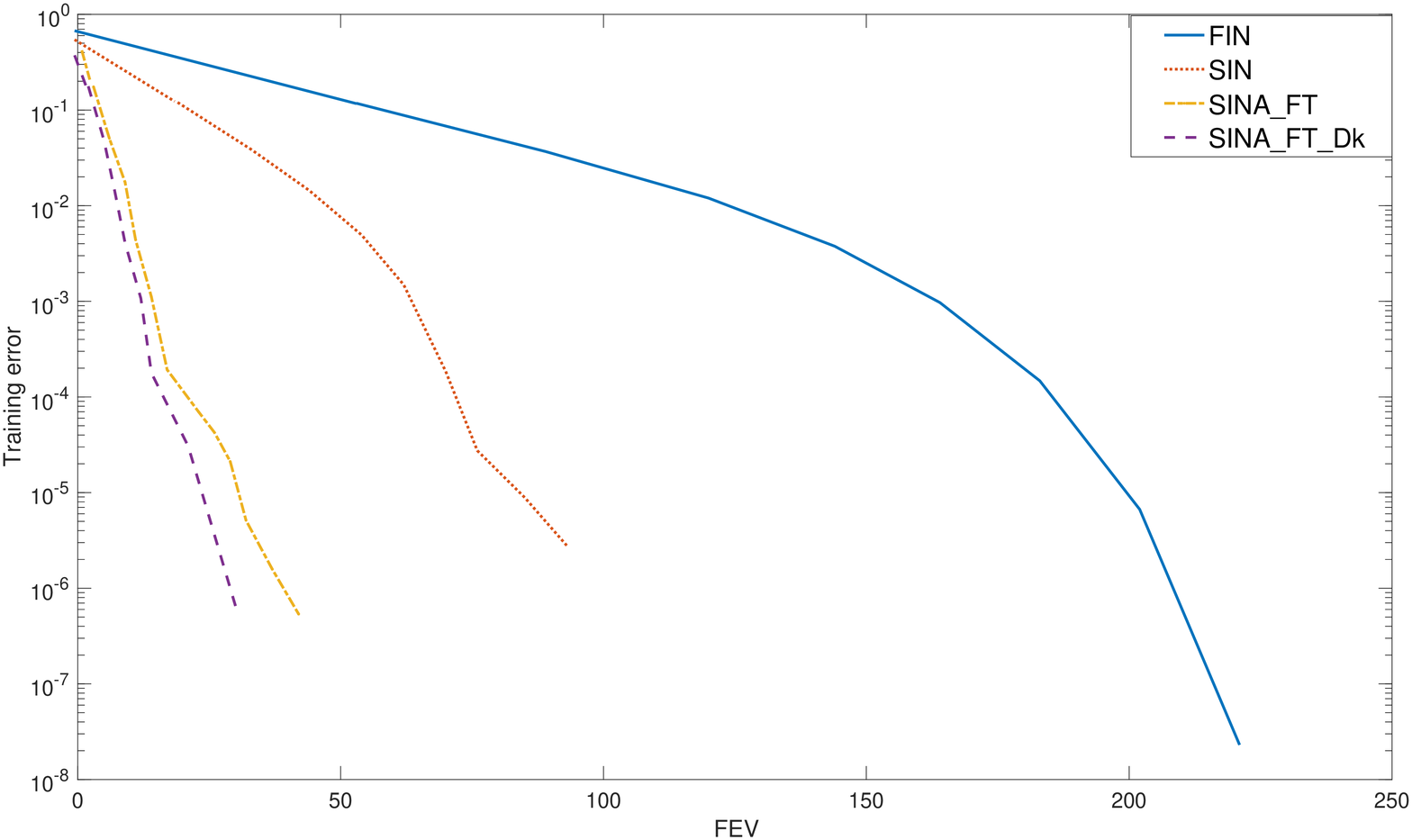}
\vspace*{-25pt}
\caption{Mushrooms dataset: training error versus iteration (left) and versus FEV (right) }
\label{training_error_mush}
\end{figure}

\begin{figure}
\hspace*{-28pt}\includegraphics[height=0.25\textheight,width=0.6\textwidth]{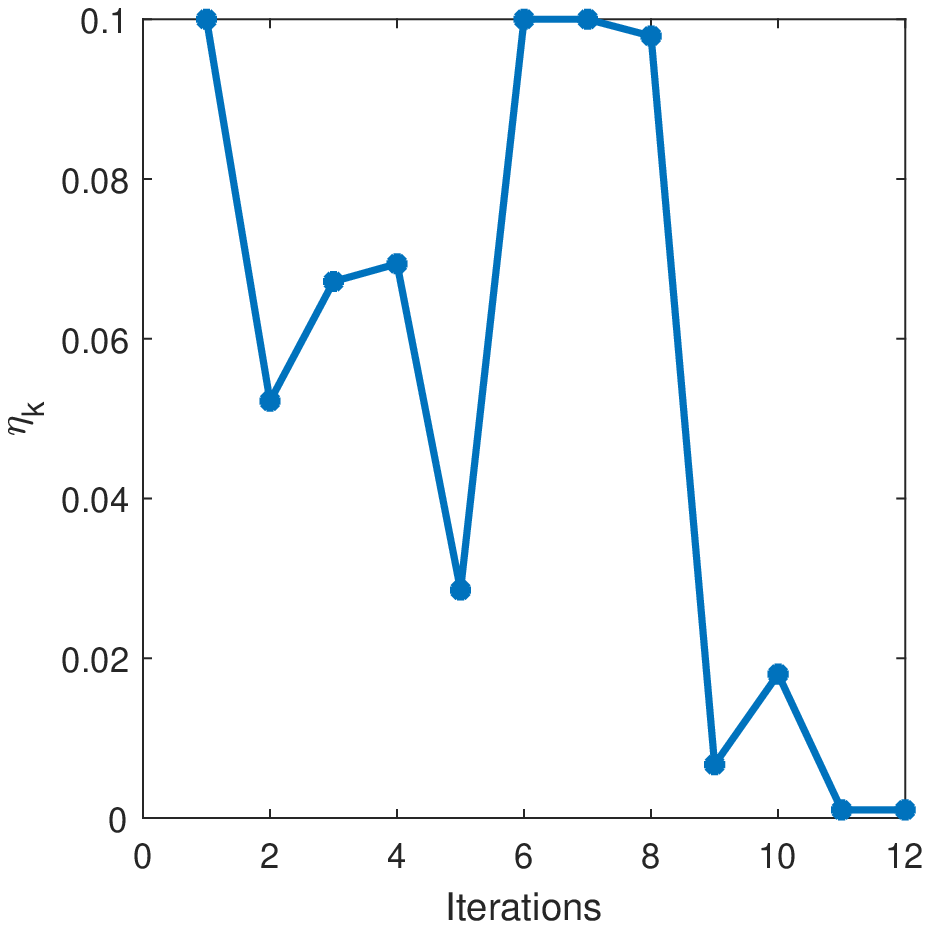}
\hspace*{-15pt}\includegraphics[height=0.25\textheight,width=0.6\textwidth]{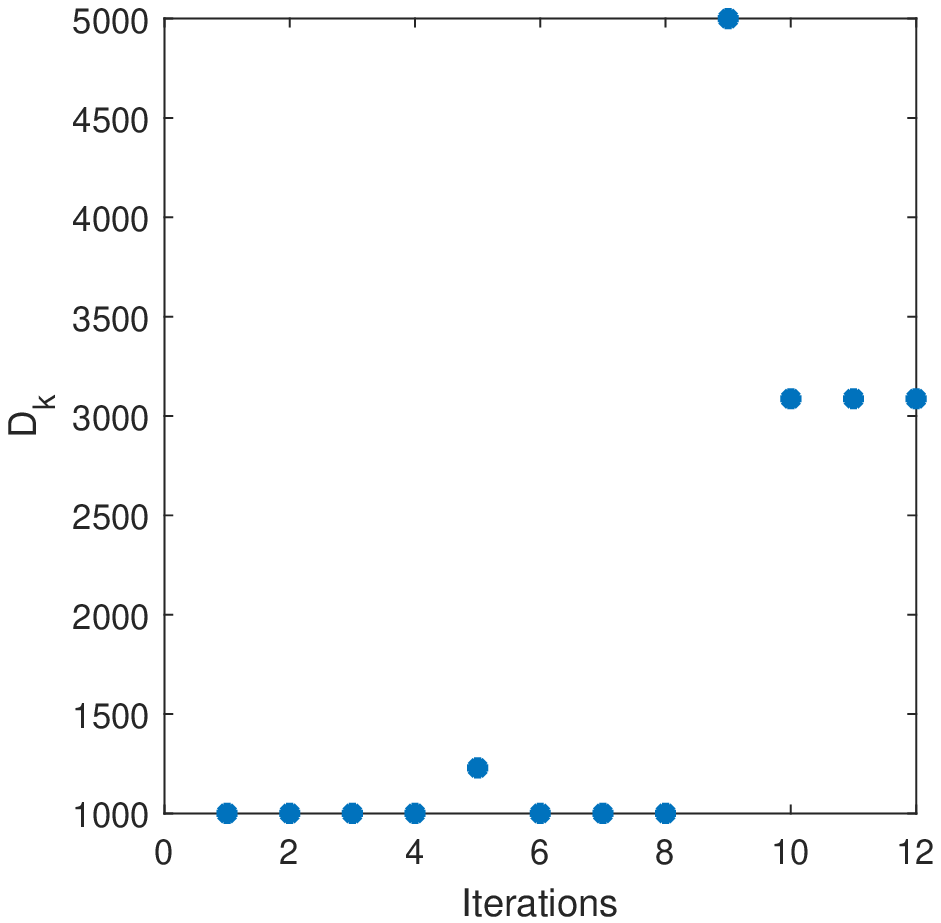}
\caption{Mushrooms dataset: {\tt SINA\_FT\_Dk}: $\eta_k$'s values  (right) an choice of $D_k$ (left) }
\label{eta_D_mush}
\end{figure}

 We begin by reporting our result with the {\tt mushrooms} dataset. We first underline that the linear algebra phase for this test is not demanding and
the average number of required  CG iteration is small, as an example it is around 10 when the adaptive choice of the forcing term is used.
In Figure \ref{training_error_mush} we plot  $f_k-f(x^*)$ (training error) versus iterations (left) and versus function evaluations FEV (right). We can observe that as expected  the full Inexact Newton (FIN) method is the fastest procedure. Moreover, the adaptive choice of the forcing terms seems to speed up the subsampled inexact procedure.
Finally, remarkably, the procedure employing both the adaptive choice of $\eta_k$ and $D_k$ seems to work quite well. Indeed it is slower than  {\tt SINA\_FT}
in the first stage of the convergence history as it uses a smaller Hessian sample set. In the last stage of the procedure it becomes faster as the sample size increases (see Figure \ref{eta_D_mush}).  On the other hand, if we look to the computational cost of the procedures we can observe that {\tt SINA\_FT\_Dk} outperforms all the procedures under comparison. We also note that  {\tt SINA\_FT}  outperforms  {\tt SIN}. Overall these results show the efficiency of the proposed adaptive strategies.

To give more insight on our adaptive choices we plot in Figure \ref{eta_D_mush}
 the values of the forcing terms $\eta_k$'s (left) and the value of $D_k$ (right) versus iterations
using the {\tt SINA\_FT\_Dk}.  We can observe that  whenever $\eta_k$ becomes smaller
$D_k$ increases and  the approximation of the hessian improves.
 Moreover, we note that the adaptive procedure allows the Hessian sample size to decrease when the model does not approximate sufficiently well the function and correspondingly the forcing term is increased.

\begin{figure}
 \hspace*{-28pt}\includegraphics[height=0.25\textheight,width=0.6\textwidth]{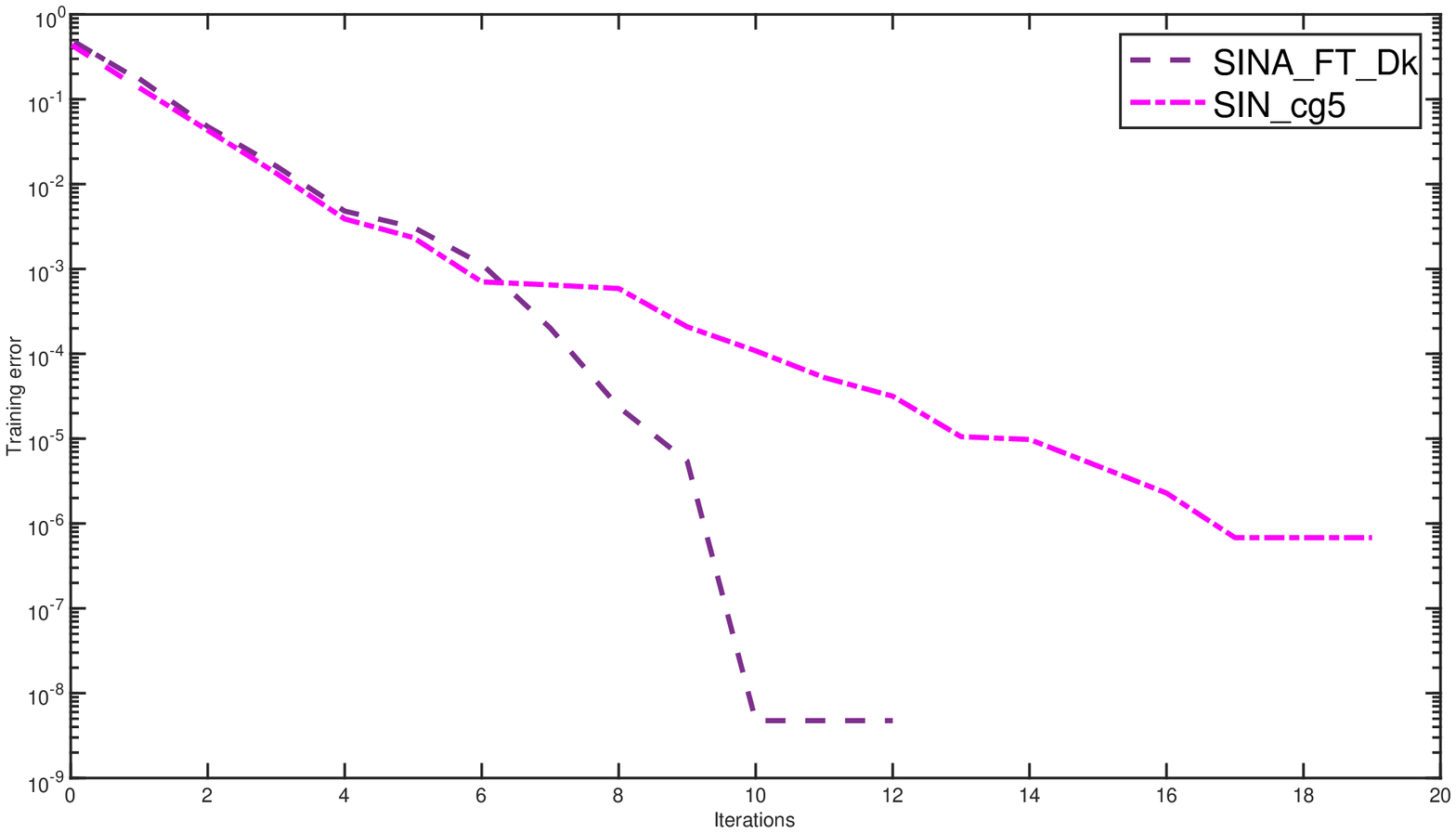}
\hspace*{-15pt} \includegraphics[height=0.25\textheight,width=0.6\textwidth]{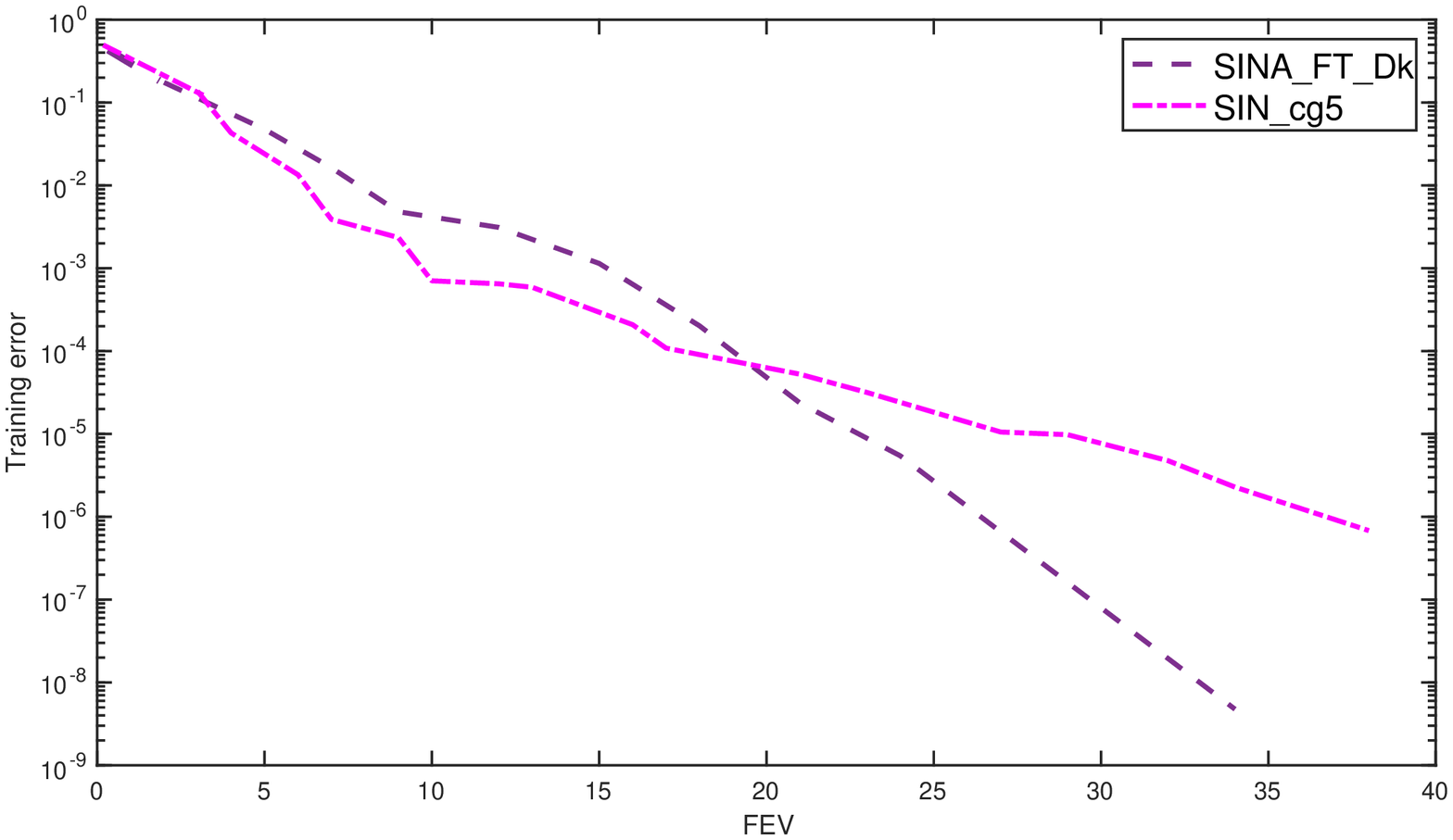}
\caption{Mushrooms dataset: training error versus iteration (left) and versus FEV (right) }
\label{mush_cg5}
\end{figure}

We also compare, in Figure \ref{mush_cg5}, our adaptive procedure {\tt SINA\_FT\_Dk} with {\tt SIN\_cg5}.
It is interesting to note that the adaptive procedure is faster than  {\tt SIN\_cg5} and greater accuracy can be reached.  In terms of FEV {\tt SIN\_cg5} is slightly better than {\tt SINA\_FT\_Dk}   till the accuracy of $10^{-4}$. If more accuracy is needed the adaptive procedure is clearly preferable.

\begin{figure}
\hspace*{-28pt}\includegraphics[height=0.25\textheight,width=0.6\textwidth]{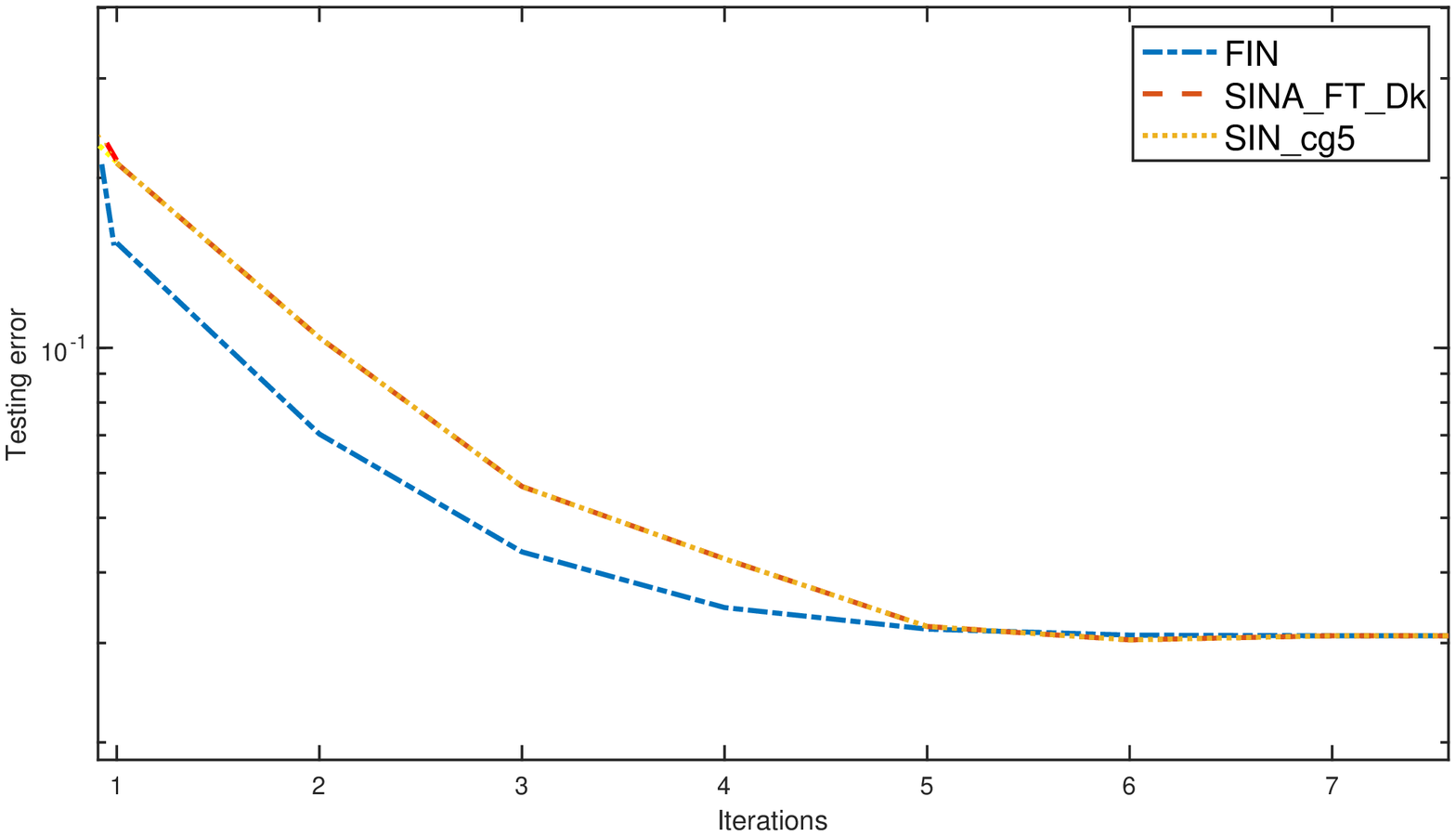}
\hspace*{-15pt}\includegraphics[height=0.25\textheight,width=0.6\textwidth]{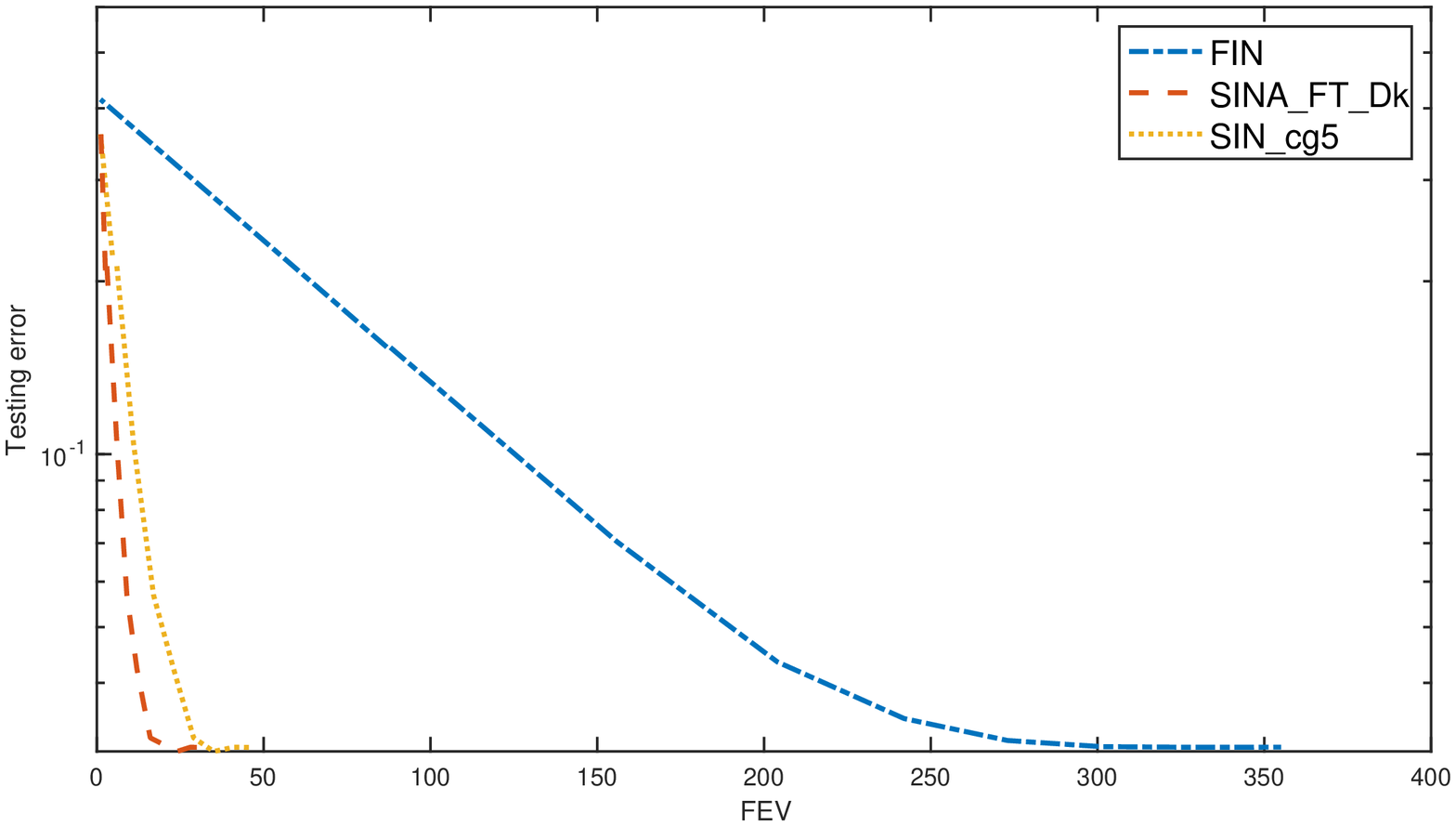}
\caption{Mushrooms dataset: testing error versus iteration (left) and versus FEV (right) }
\label{testing_error_mushrooms}
\end{figure}

Finally, in Figure \ref{testing_error_mushrooms}  we compare the behaviour of  testing errors versus iterations (left) and   FEVs (right) along the sequences generated by  full Newton method {\tt FIN},
the adaptive procedure {\tt SINA\_FT\_Dk} and  {\tt SIN\_cg5} method.
We evaluated the testing error as follows.
 Let $x^k$ be the approximation computed  by  a method using the data in the training set. Then, $x^k$  is used to classify the samples in the testing set made up of $\bar N=3124$ instances $z_i$, $i=1,\ldots,\bar N$ and corresponding $b_i$'s. The classification error at iteration $k$  is defined as $\frac{1}{\bar N}\sum_{i=1}^{\bar N}\log(1+\exp(-b_i z_i^T x^k))$.
  We can observe that the three methods reached  almost the same testing errors, the subsampled approaces require a lower computational cost and again the adaptive procedure outperforms {\tt SIN\_cg5}.
   We also observe that the testing error quickly steadies. This enlightens that the stopping tolerance we used is tight enough and going on with the iterations 
  would not decrease further the testing error.

\begin{figure}
\hspace*{-28pt} \includegraphics[height=0.25\textheight,width=0.6\textwidth]{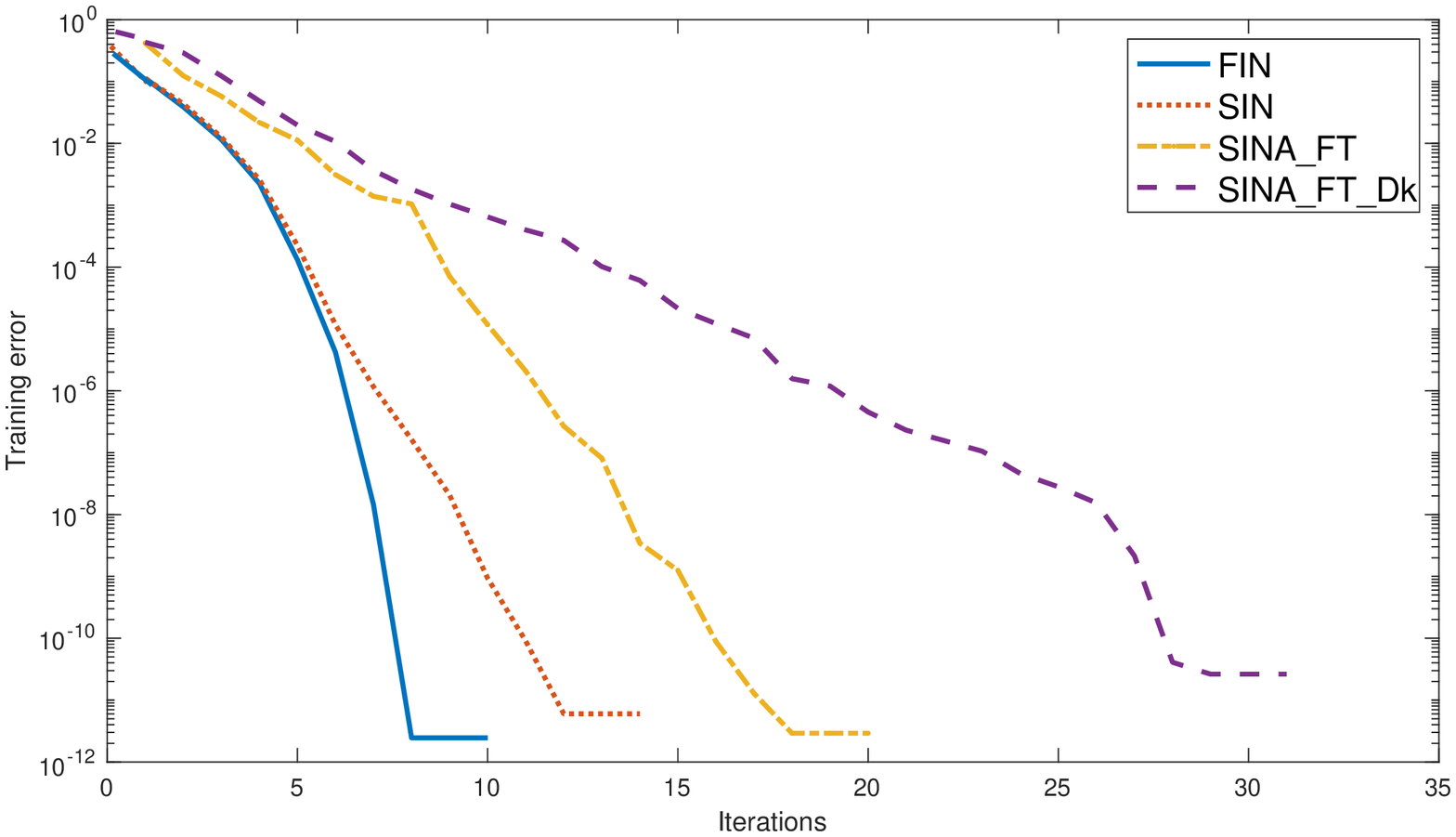}
\hspace*{-15pt} \includegraphics[height=0.25\textheight,width=0.6\textwidth]{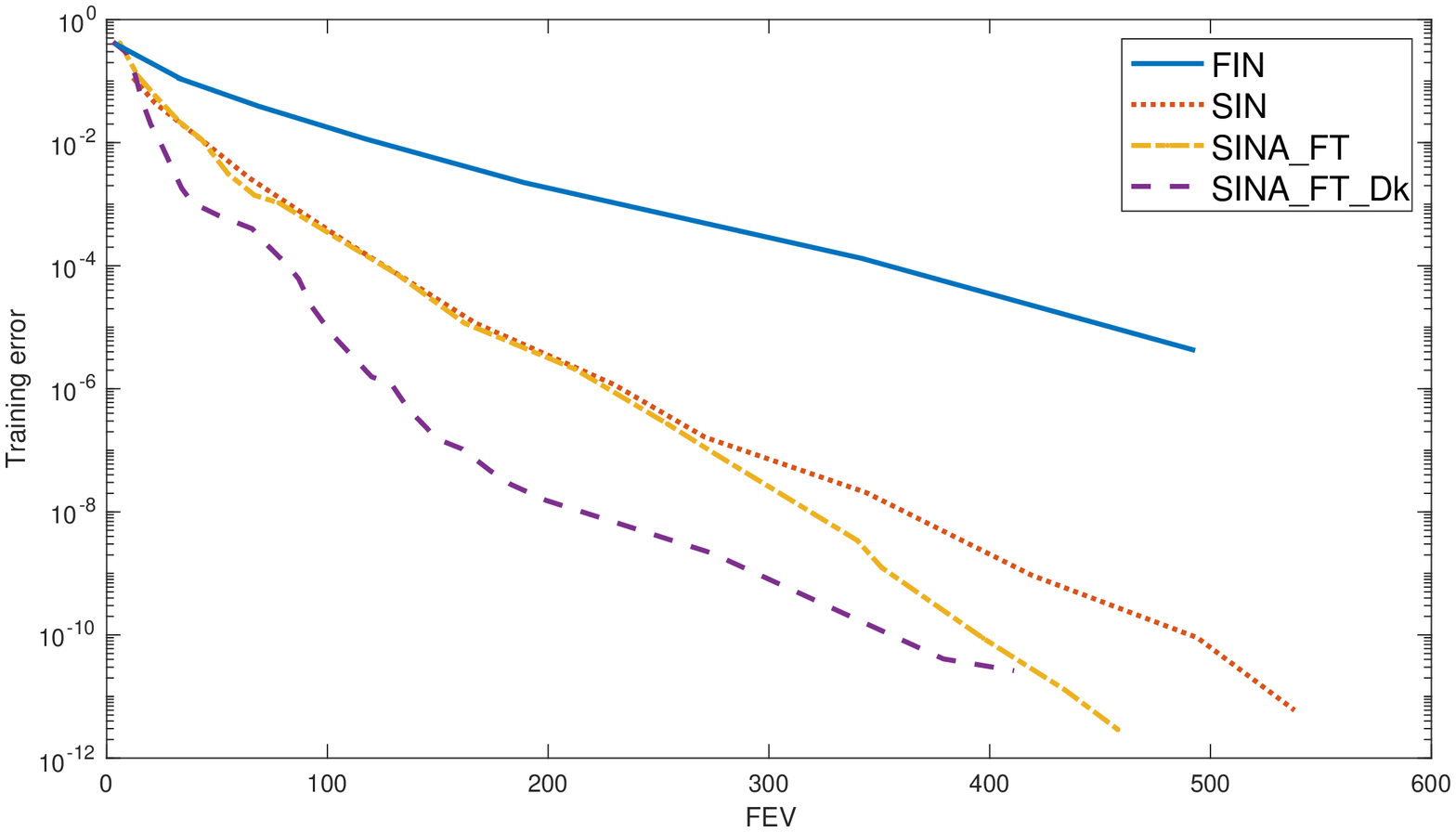}
\caption{CINA0 dataset: training error versus iteration (left) and versus FEV (right) }
\label{training_error_CINA0}
\end{figure}

\begin{figure}
\hspace*{-28pt} \includegraphics[height=0.25\textheight,width=0.6\textwidth]{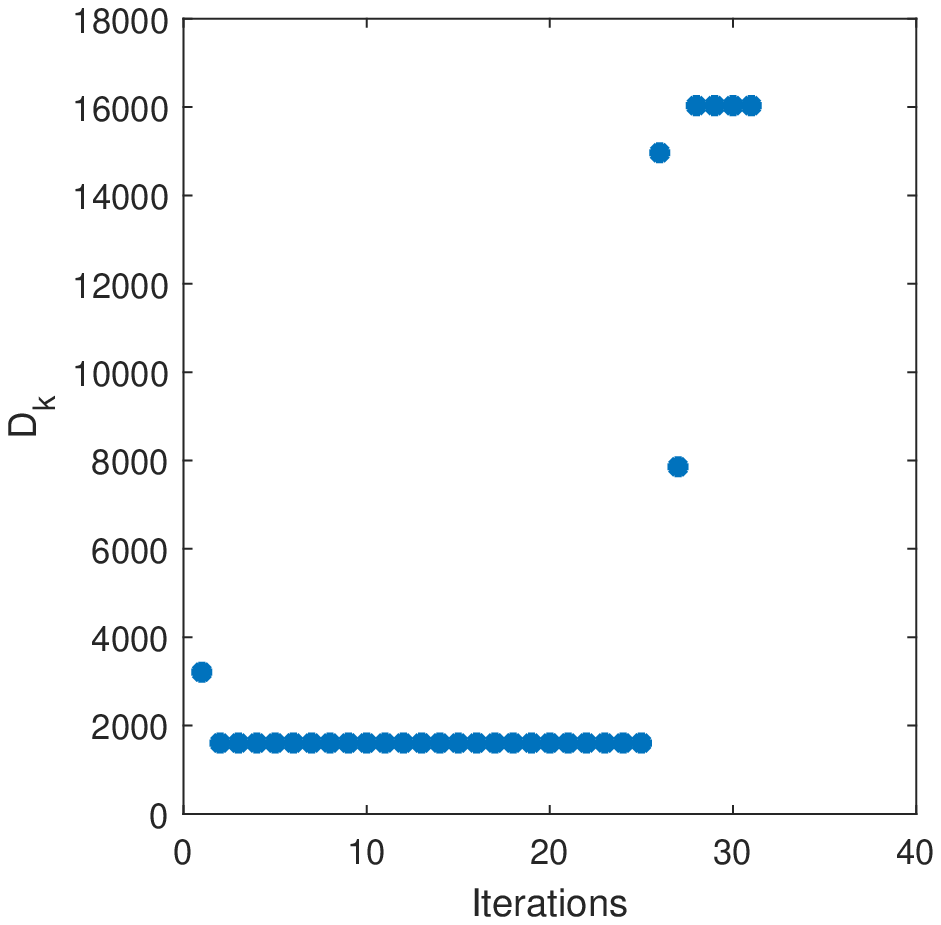}
\hspace*{-15pt}\includegraphics[height=0.25\textheight,width=0.6\textwidth]{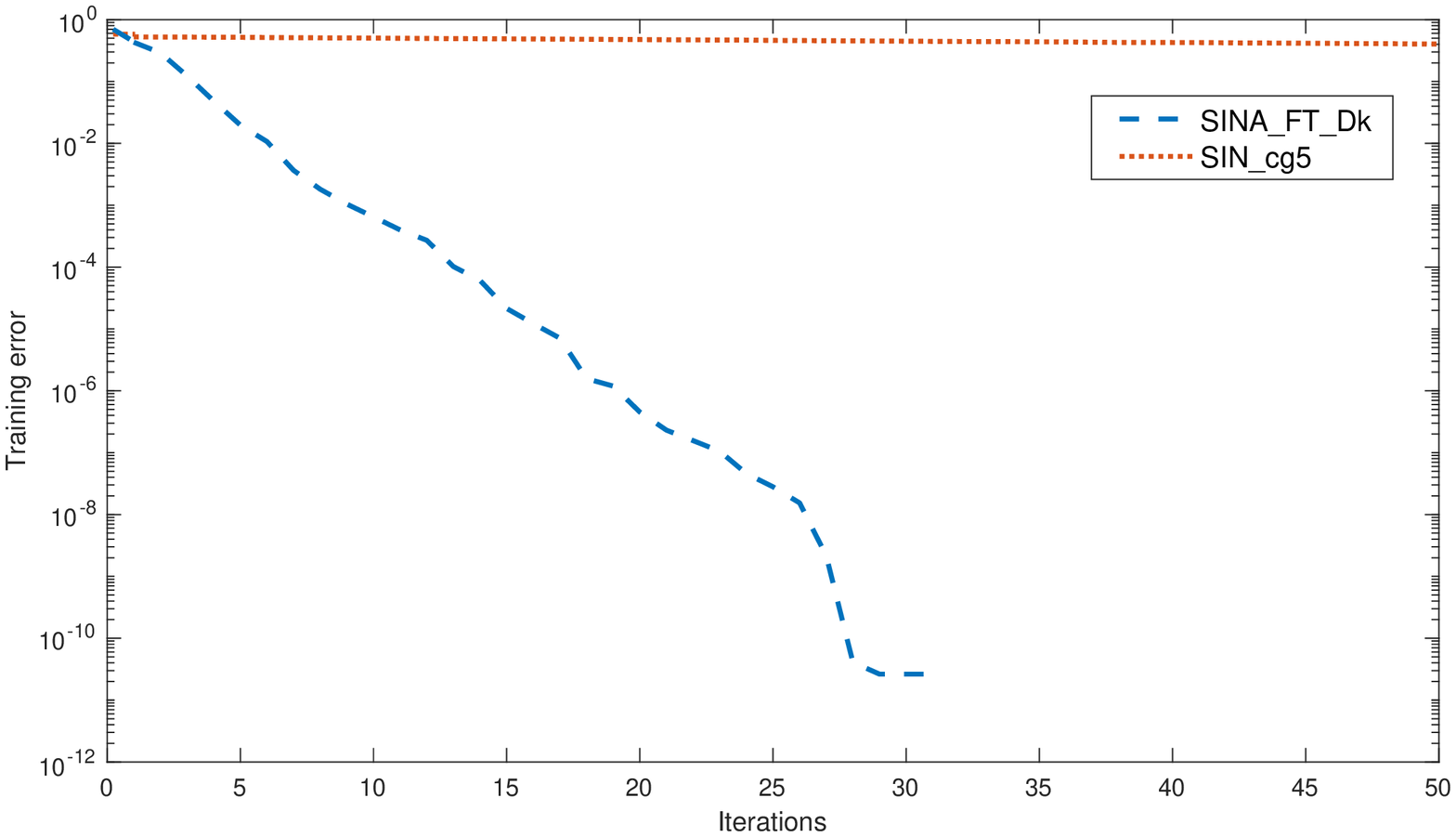}
\caption{CINA0 : $D_k$ values in {\tt SINA\_FT\_DK} (left),  training error (right) }
\label{CINA0_DK}
\end{figure}

Let us consider the CINA0 dataset.
In Figure \ref{training_error_CINA0} we plot  the training error versus iterations (left) and versus F-evaluations (right). The full Inexact Newton method {\tt FIN}  is the fastest one, as expected. On the other hand, also in this case, it is  the most expensive. The behaviour of the subsampled procedures is that desiderable, with the adaptive procedures that outperforms the {\tt SINA\_FT} in terms of FEV. Note that {\tt SINA\_FT\_Dk} is slower than the other two subsampled procedures, but less costly,  as expected. In fact this is a problem where the linear algebra phase is more demanding, the average number of CG iteration using
{\tt SINA\_FT}  is about 70, and small values of $D_k$ are used to limit the overall computational cost (see Figure \ref{CINA0_DK}, left). Since the convergence of CG is slow {\tt SIN\_cg5} is not able to converge with a reasonable rate and it can only provide a very rough accuracy (see Figure   \ref{CINA0_DK}, right)

\begin{figure}
\hspace*{-28pt} \includegraphics[height=0.25\textheight,width=0.6\textwidth]{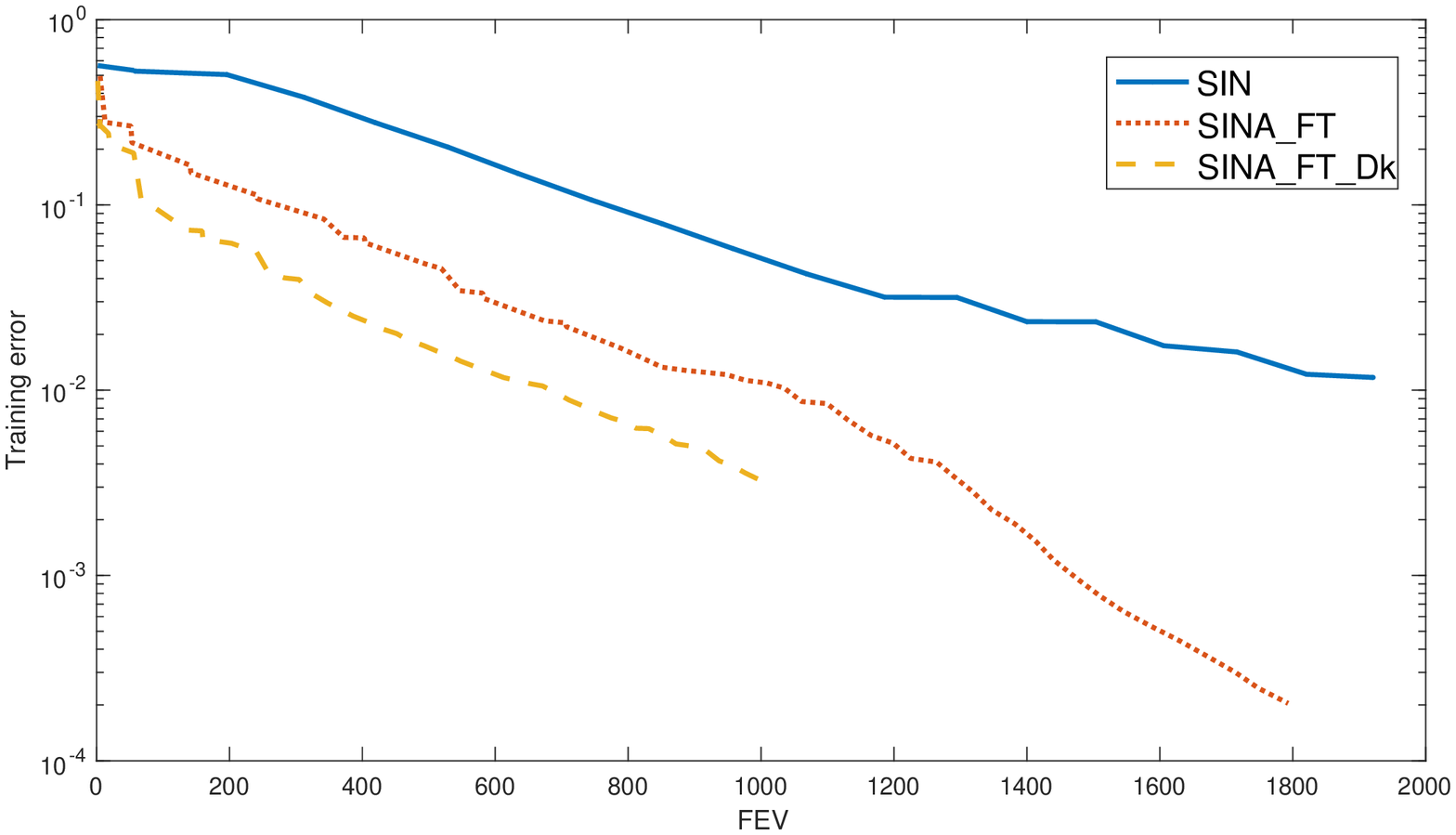}
\hspace*{-15pt}  \includegraphics[height=0.25\textheight,width=0.6\textwidth]{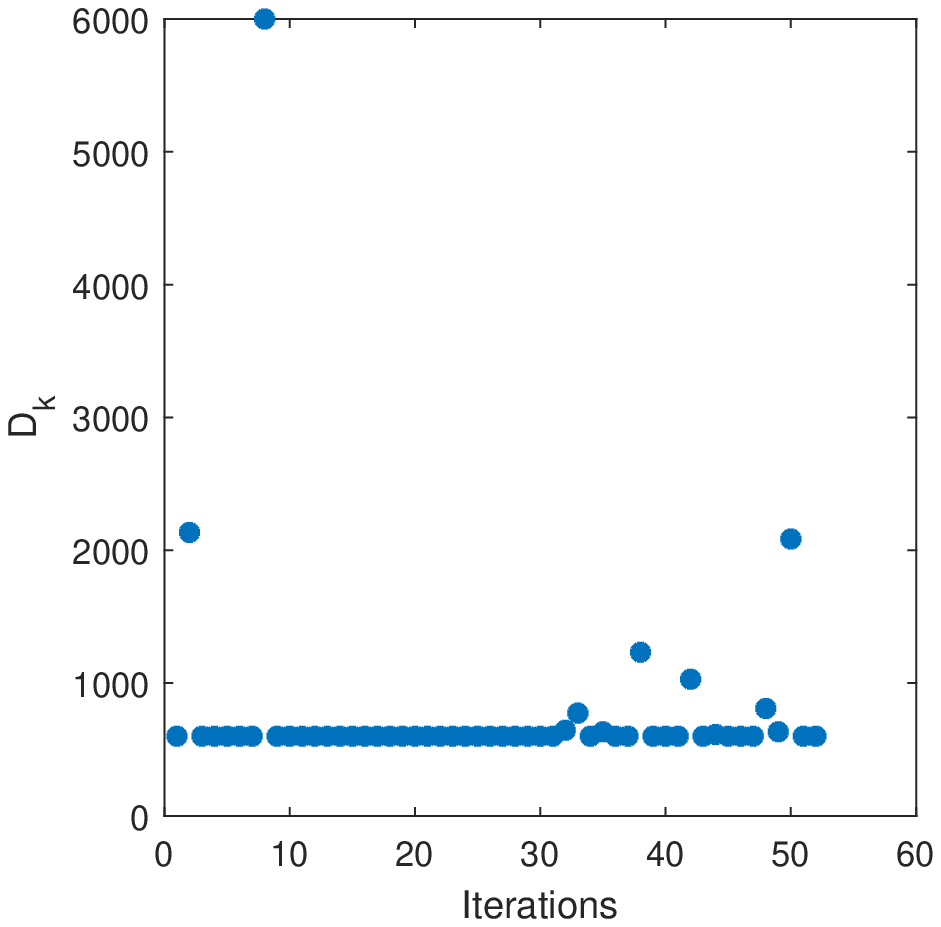}
\caption{Gisette dataset: training error versus FEV (left), $D_k$'s versus Newton iteration (right), $\eta_k$'s versus Newton iterations (bottom-left) and cg-iterations versus Newton iterations (bottom-right) }
\label{training_error_Gisette_01}
\end{figure}

\begin{figure}
 \hspace*{-28pt}\includegraphics[height=0.25\textheight,width=0.6\textwidth]{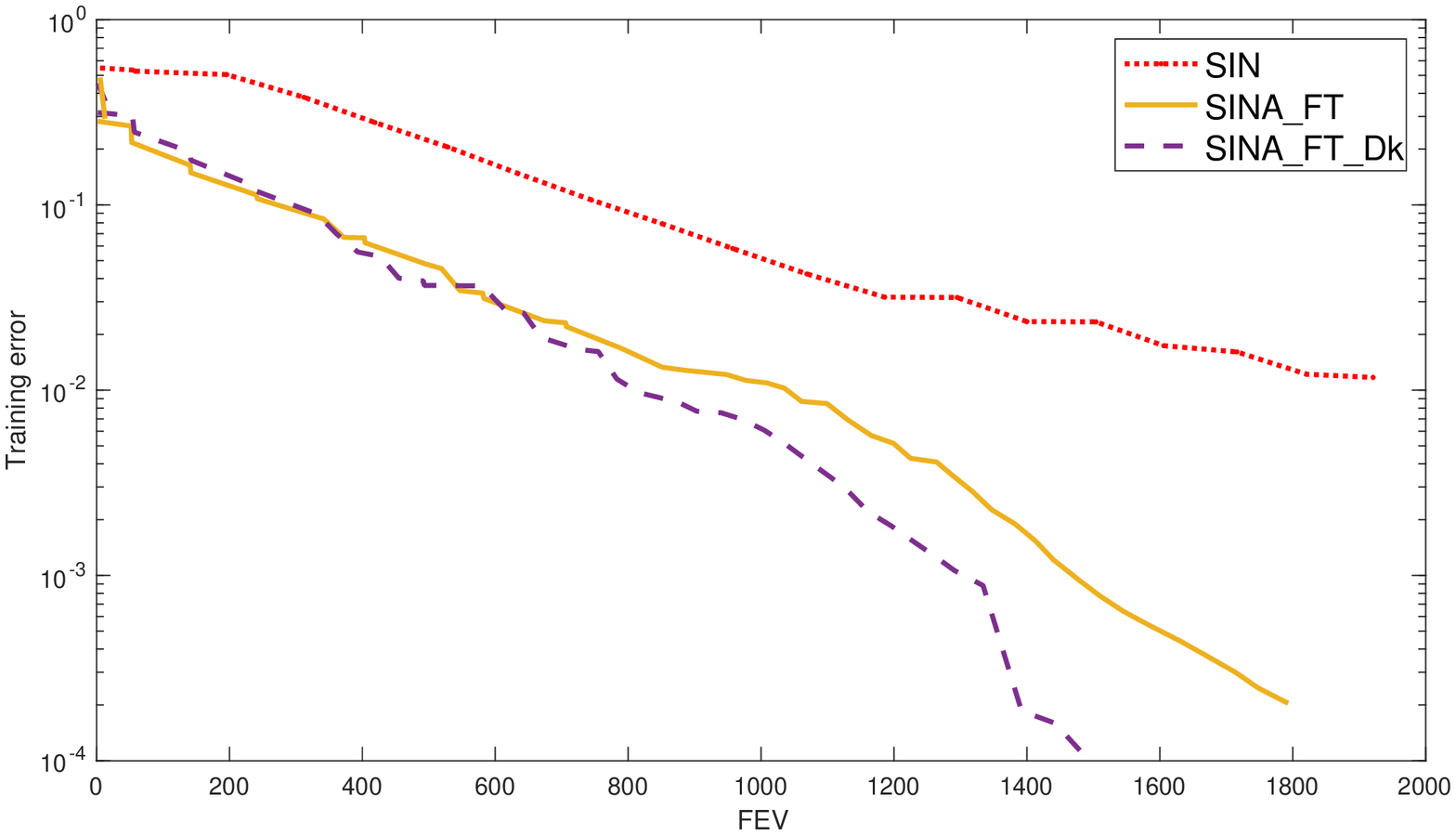}
 \hspace*{-15pt}\includegraphics[height=0.25\textheight,width=0.6\textwidth]{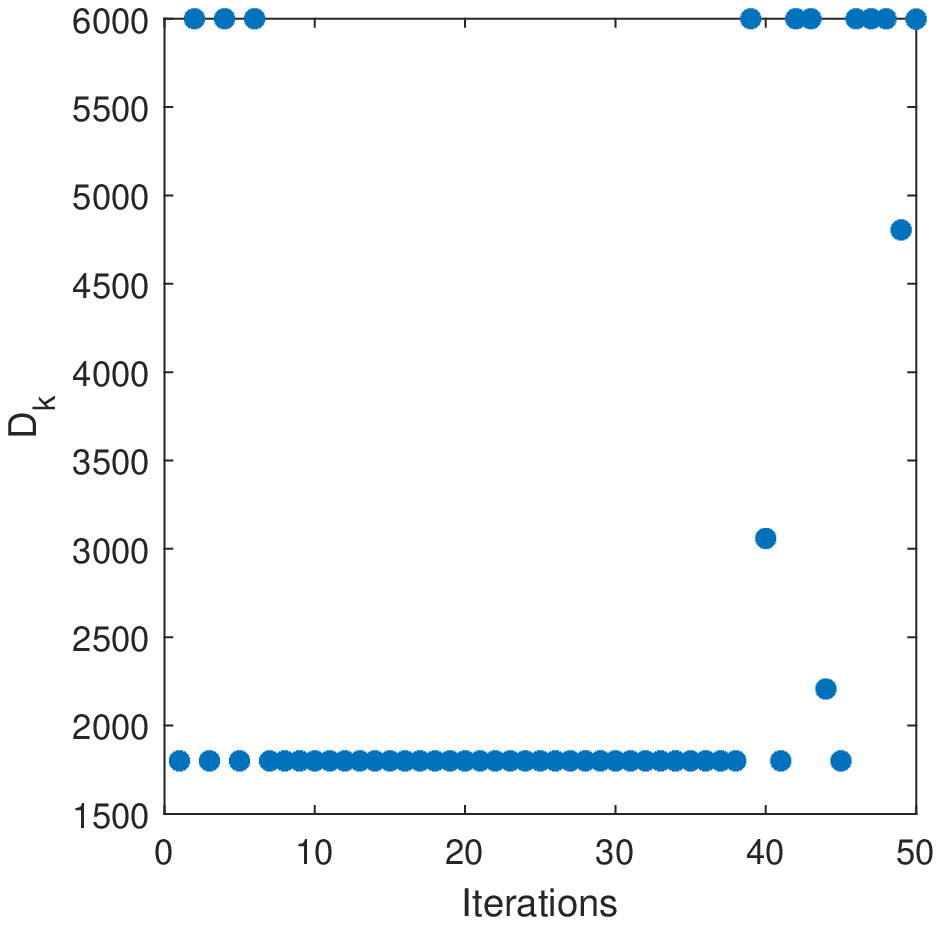}
\hspace*{-28pt} \includegraphics[height=0.25\textheight,width=0.6\textwidth]{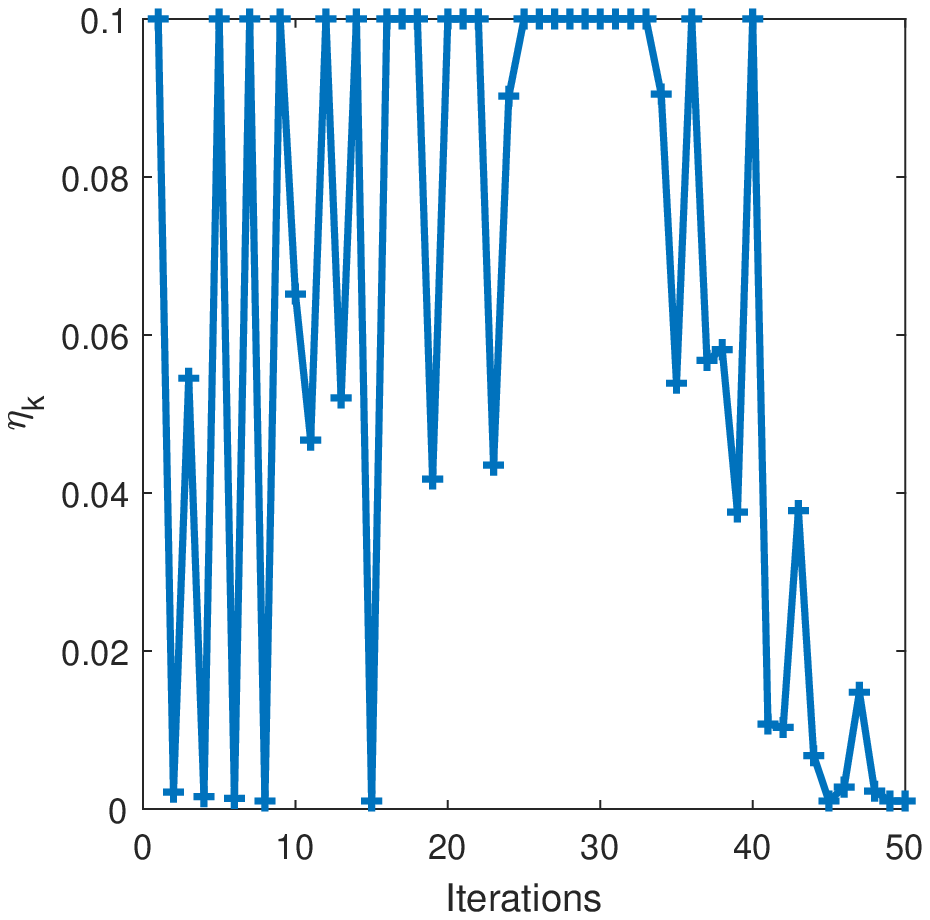}
\hspace*{-15pt} \includegraphics[height=0.25\textheight,width=0.6\textwidth]{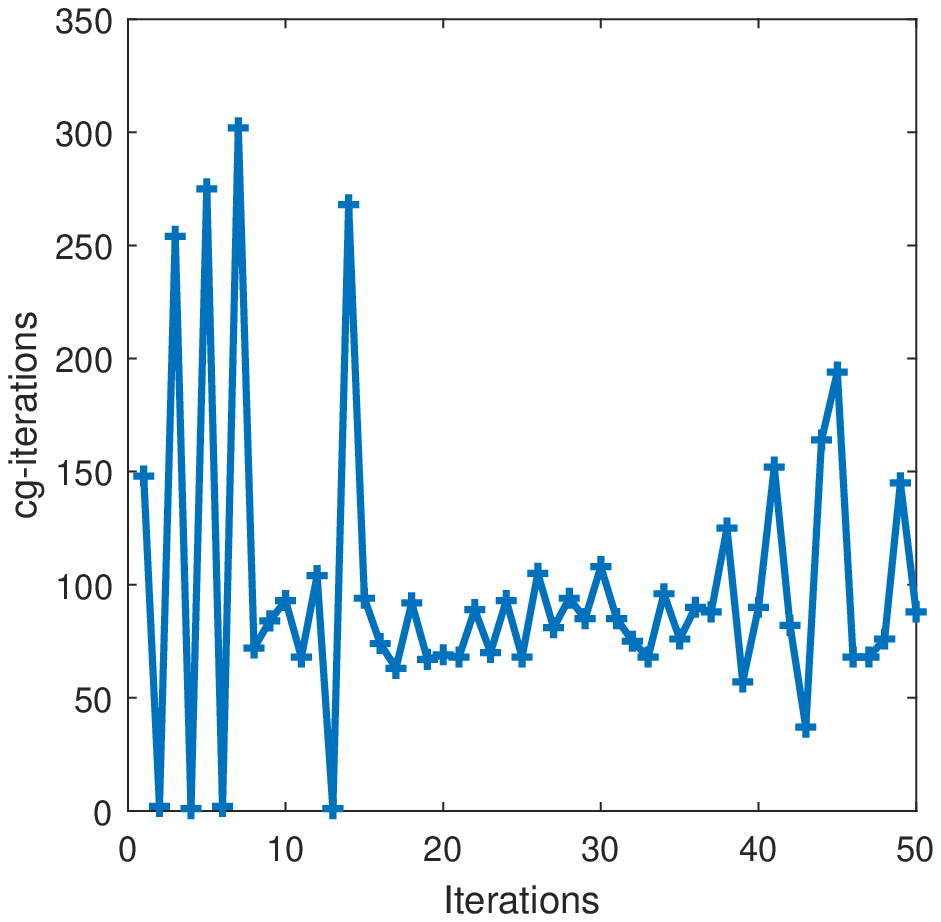}
\caption{Gisette dataset ($D_0=0.3N$ in  {\tt SINA\_FT\_DK}): training error versus FEV (top-left), $D_k$'s versus Newton iteration (top-right), $\eta_k$'s versus Newton iterations (bottom-left) and cg-iterations versus Newton iterations (bottom-right) }
\label{training_error_Gisette}
\end{figure}

We finally show the results obtained with the Gisette dataset, where $n$ is larger than in the previous tests.  In Figure  \ref{training_error_Gisette_01}, we compare  {\tt SINA\_FT\_DK},  {\tt SIN} and  {\tt SIN\_FT} in terms of   function evaluations FEV (left) and we also report the behaviour of the sample size along the iterations (right). We can observe that     {\tt SIN\_FT} and  {\tt SINA\_FT\_DK}  are  less expensive and more accurate
than {\tt SIN}, however within 50 nonlinear iterations {\tt SINA\_FT\_DK}   is not able to produce an approximation as accurate as that provided by   {\tt SIN\_FT}.
In fact, the convergence is slow due to the small size of the Hessian subsample set.  However, note that if a training error of the order of $10^{-2}$  is enough,  {\tt SINA\_FT\_DK} is the method of choice as it is less expensive. In case a greater accuracy is needed $D_0=0.3N$ should be used in  {\tt SINA\_FT\_DK}.  Figure  \ref{training_error_Gisette} refers to this choice of $D_0$.   In this case {\tt SINA\_FT\_DK} outperforms both  {\tt SIN} and   {\tt SIN\_FT}.
 We also report the values of $D_k$, $\eta_k$ and the number of cg-iterations at each Newton iteration.  Plots show that we use small values of $D_k$ and large values of $\eta_k$ till the last stage of convergence. Since we are solving linear systems with a rough accuracy, the number of cg-iterations is reasonable, except for a few occurences,  considering that $n=5000$ and we are not employing a preconditioner.
\vskip 10 pt 
\noindent
{\bf Acknowledgements}\\
The authors are grateful to the anonymous referees for their suggestions, which lead to significant
improvement of the manuscript.

\section{Appendix}

\noindent
{\bf Proofs of Lemmas  \ref{Lemma4},  \ref{Lemma5} and \ref{dkbounds}}

\noindent
{\bf  Proof of Lemma \ref{Lemma4}.}
Using the Mean Value Theorem  we get the following inequality for the full Hessian.
$$\|\nabla f_{\cal N}(x^{k}+s^k)-\nabla f_{\cal N}(x^k) - \nabla^2 f_{\cal N}(x^k)s^k\|\leq \frac{L}{2}\|s^k\|^2.$$
Now,
\begin{eqnarray*}
\|\nabla f_{\cal N}(x^{k}+s^k)-\nabla f_{\cal N}(x^k) - \nabla^2 f_{\cal D}(x^k)s^k\| & \leq  & \|\nabla f_{\cal N}(x^{k}+s^k)-\nabla f_{\cal N}(x^k) - \nabla^2 f_{\cal N}(x^k)s^k\| \\
 & + & \| \nabla^2 f_{{\cal D}_k}(x^k)-\nabla^2 f_{\cal N}(x^k)\| \|s^k\| \\
 & \leq & \|s^k\|  (h(D_k,x^k) + \frac{1}{2} L \|s^k\|).
\end{eqnarray*}
Then,
\begin{eqnarray*}
\|\nabla f_{\cal N}(x^{k}+s^k) \| & \leq  & \|\nabla f_{\cal N}(x^k) + \nabla^2 f_{{\cal D}_k}(x^k)s^k\| + \|\nabla f_{\cal N}(x^{k}+s^k) - \nabla f_{\cal N}(x^k) - \nabla^2 f_{{\cal D}_k}(x^k) s^k\| \\
& \leq & \eta \|\nabla f_{\cal N}(x^k)\| +\|s^k\|(\frac{1}{2}L \|s^k\| + h(D_k,x^k))
\end{eqnarray*}
and, by Lemma \ref{Lemma3a},
\begin{eqnarray*}
\|\nabla f_{\cal N}(x^{k}+s^k) \| & \leq & \eta \|\nabla f_{\cal N}(x^k)\| +
\frac{1}{\lambda_1} \|\nabla f_{\cal N}(x^k)\|(\frac{1}{2\lambda_1}  L \|\nabla f_{\cal N}(x^k)\| + h(D_k,x^k))
\end{eqnarray*}
$ \Box $

\noindent
{\bf Proof of Lemma  \ref{Lemma5}.}  By Lemma \ref{Lemma3a} and (\ref{Lemma2}) we have
\begin{eqnarray*}
\|x^k+s^k - x^*\| & \leq & \|x^k - x^*\| + \|s^k\| \leq \|x^k-x^*\| + \lambda_1^{-1} \|\nabla f_{\cal N}(x^k)\| \\
& \leq & \|x^k-x^*\| + \lambda_1^{-1}  \lambda_{n} \|x^k-x^*\| \leq (1+\lambda_1^{-1} \lambda_{n})\|x^k-x^*\| \leq \delta^* .
\end{eqnarray*} $ \Box $

\noindent
{\bf Proof of Lemma \ref{dkbounds}.}
Let us define $Y_i(x)=(H_i(x)-\nabla^2 f_{\cal N}(x))/D$ where $H_i(x)$ is a randomly chosen Hessian.   Then,
$$\sum_{i=1}^{D} Y_i(x)=\nabla^2 f_{\cal D}(x)-\nabla^2 f_{\cal N}(x).$$
 Also, notice that $E(Y_i(x))=0$  and that the Weyl's inequality  yields $\lambda_{max}(Y_i(x))\leq \lambda_n/D$,  $\lambda_{min}(Y_i(x))\geq -\lambda_n/D$. 
 Then, $\|Y_i(x)\|\le \lambda_n/D$ and

$$
 \tilde{\sigma}^2 (x) := \|\sum_{i=1}^{D} E(Y^2_i(x))\| \leq   \sum_{i=1}^{D} E(\|Y_i(x)^2\|) \leq \frac{\lambda_n^2}{D}.
$$
 Then,  using  the Bernstein's inequality (see \cite[Theorem 1.6]{tropp12} )  we derive 
 $$
P(\|\nabla^2 f_{\cal D}(x)-\nabla^2 f_{\cal N}(x)\| \geq \gamma)\leq 
2 n \; \exp \left(-\frac{ \gamma^2/2}{\lambda_n^2/D+( \lambda_n/D)\gamma/3}\right).
$$
This yields the bound.
$\Box$

\end{document}